\documentclass{amsart}

\usepackage{amsmath,amsthm,amssymb}

\usepackage{graphicx}

\newcommand{\alb}{\mathsf{X}}

\newcommand{\arr}{\longrightarrow}

\newcommand{\bim}{\mathfrak{M}}
\newcommand{\C}{\mathbb{C}}
\newcommand{\CS}{{\widehat\C}}
\newcommand{\CP}{\mathbb{C}\mathbb{P}}
\newcommand{\dom}{\mathop{\mathrm{Dom}}}

\newcommand{\F}{\mathcal{F}}

\newcommand{\img}[1]{\mathop{\mathrm{IMG}}\left(#1\right)}
\newcommand{\lims}[1][G]{{\mathcal{J}_{#1}}}
\newcommand{\limg}[1][G]{\mathcal{X}_{#1}}
\newcommand{\til}{\mathcal{T}}
\newcommand{\M}{\mathcal{M}}
\newcommand{\mapdown}[1]{\Big\downarrow\rlap{$\vcenter{\hbox{$\scriptstyle#1$}}$}}
\newcommand{\N}{\mathbb{N}}
\newcommand{\nuke}{\mathcal{N}}
\newcommand{\R}{\mathbb{R}}
\newcommand{\si}{\mathsf{s}}
\newcommand{\symm}[1][\alb]{\mathfrak{S}\left(#1\right)}
\newcommand{\wt}{\widetilde}
\newcommand{\X}{\mathcal{X}}
\newcommand{\xmo}{\alb^{-\omega}}

\newcommand{\xs}{\alb^*}
\newcommand{\Z}{\mathbb{Z}}

\newcommand{\hide}[1]{}

\newtheorem{theorem}{Theorem}[section]
\newtheorem{proposition}[theorem]{Proposition}
\newtheorem{corollary}[theorem]{Corollary}
\newtheorem{lemma}[theorem]{Lemma}

\theoremstyle{definition}
\newtheorem{defi}{Definition}
\newtheorem{examp}{Example}
\newtheorem*{remark}{Remark}

\title{Combinatorial models of expanding dynamical systems}
\author{Volodymyr Nekrashevych}
\thanks{This material is based upon work supported by the
National Science Foundation under Grants DMS0605019 and
DMS0800085.}

\begin{document}
\maketitle \tableofcontents

\begin{abstract}
We define iterated monodromy groups of more general structures
than partial self-covering. This generalization makes it possible
to define a natural notion of a combinatorial model of an
expanding dynamical system. We prove that a naturally defined
``Julia set'' of the generalized dynamical systems depends only on
the associated iterated monodromy group. We show then that the
Julia set of every expanding dynamical system is an inverse limit
of simplicial complexes constructed by inductive cut-and-paste
rules.
\end{abstract}

\section{Introduction}

According to a well known principle, expanding (and more
generally, hyperbolic) dynamical systems have combinatorial nature
and are determined by a finite amount of data. For instance, they
are \emph{finitely presented}, see~\cite{fried}. This principle
can be also formulated in a form of structural stability or
rigidity theorems: two hyperbolic dynamical systems that are
topologically or homotopically close to each other are conjugate.

The aim of our paper is to describe, by proving the corresponding
rigidity theorem, a complete algebraic invariant of expanding
dynamical systems. We translate then this algebraic invariant into
a more geometric language of polyhedral models of dynamical
systems (and their Julia sets). These models give a representation
of the Julia set of the dynamical system as an inverse limit of
simplicial complexes that are constructed using simple
cut-and-paste rules, similar to subdivision rules in
one-dimensional complex dynamics. We illustrate our techniques, in
particular, by constructing combinatorial models of the Julia sets
of multi-dimensional dynamical systems.

We define the algebraic invariant (called the \emph{iterated
monodromy group}) in a general setting of a multi-valued partially
defined dynamical system. Namely, a \emph{topological automaton}
is a pair of maps $f:\M_1\arr\M$, $\iota:\M_1\arr\M$ between two
topological spaces (or orbispaces), such that $f$ is a finite
covering map. If $\iota$ is a homeomorphism, then we can identify
$\M_1$ and $\M$ by $\iota$, thus getting a dynamical system
$f:\M\arr\M$. If $\iota$ is an embedding, then $f:\M_1\arr\M$ is a
partial self-covering.

Iterated monodromy groups were originally defined for partial
self-coverings only (see~\cite{bgn,nek:book}). However, the fact
that $\iota$ is an embedding is not used neither in the
construction nor in the main results of~\cite{nek:book}. Moreover,
iterated monodromy groups of partial self-coverings of orbispaces
are defined in~\cite{nek:book} essentially in the setting of
topological automata.

Topological automata were studied (under different names) by
T.~Katsura in~\cite{katsura:generalizing} in relation with
$C^*$-algebras, and by Y.~Ishii and J.~Smillie~\cite{ishiismillie}
in relation with homotopical rigidity of hyperbolic dynamical
systems. The last article was one of inspirations of our paper.

A topological automaton $f, \iota:\M_1\arr\M$ can be naturally
iterated. The covering $f$ and the map $\iota$ induce a covering
$f_1:\M_2\arr\M_1$ and a map $\iota_1:\M_2\arr\M_1$ defined by the
pull-back diagram
\[\begin{array}{ccc}\M_2 & \stackrel{\iota_1}{\arr} & \M_1\\
\mapdown{f_1} &  & \mapdown{f}\\
\M_1 & \stackrel{\iota}{\arr} & \M\end{array}\] We define then,
inductively, coverings $f_n:\M_{n+1}\arr\M_n$ and maps
$\iota_n:\M_{n+1}\arr\M_n$. The $n$th iteration of the pair $f,
\iota:\M_1\arr\M$ is then the pair \[f\circ f_1\circ\cdots\circ
f_{n-1},\quad
\iota\circ\iota_1\circ\cdots\circ\iota_{n-1}:\M_n\arr\M.\] If
$\iota$ is an embedding, then the spaces $\M_n$ are the domains of
the iterations $f^n$ of the partial map $f:\M_1\arr\M$.

The iterated monodromy of a topological automaton is defined in
Section~\ref{s:img}. Rather than to give the definition here, we
define an equivalent notion of the associated virtual endomorphism
of the fundamental group. Suppose that the space $\M$ is path
connected and locally path connected. Since $f:\M_1\arr\M$ is a
finite covering map, the induced map
$f_*:\pi_1(\M_1)\arr\pi_1(\M)$ is an embedding, and
$f_*(\pi_1(\M_1))$ has finite index in $\pi_1(\M)$. The
\emph{virtual endomorphism} associated with the topological
automaton $f, \iota:\M_1\arr\M$ is the homomorphism $\iota_*\circ
f_*^{-1}$ from the subgroup $f_*(\pi_1(\M_1))\le\pi_1(\M)$ to
$\pi_1(\M)$. It is well defined up to inner automorphisms of
$\pi_1(\M)$.

If two topological automata $f', \iota':\M_1'\arr\M'$ and $f'',
\iota'':\M_1''\arr\M''$ have the same associated virtual
endomorphisms $\phi'$ and $\phi''$ (i.e., if there exists an
isomorphism $\alpha:\pi_1(\M')\arr\pi_1(\M'')$ such that
$\alpha\circ\phi'$ is equal to $\phi''\circ\alpha$ modulo inner
automorphisms), then the topological automata are called
\emph{combinatorially equivalent}. More generally, the automata
are combinatorially equivalent, if we can make the associated
virtual endomorphisms the same by taking quotients of the
fundamental groups by normal subgroups invariant under the action
of the virtual endomorphisms. A more precise definition is given
in Subsection~\ref{ss:combequiv}.

A topological automaton $f, \iota:\M_1\arr\M$ is called
\emph{contracting} if $\M$ and $\M_1$ are compact length spaces
(e.g., Riemannian manifolds, or simplicial complexes with
Riemannian structure on simplices), $f$ is a local isometry, and
$\iota$ is contracting. If $\iota$ is a homeomorphism or an
embedding, then it could be more natural to restrict the length
structure of $\M$ onto $\M_1$. In this setting an equivalent
condition is that $f$ is \emph{expanding}.

Our first main result is the following rigidity theorem (see
Theorems~\ref{th:approximationlimg} and~\ref{th:projlim}).

\begin{theorem}
\label{th:main1} Let $\F=(\M, \M_1, f, \iota)$ be a contracting
topological automaton with semi-locally simply connected space
$\M$. Denote by $\lim_\iota\F$ the inverse limit of the sequence
of spaces and maps
\[\M\stackrel{\iota}{\longleftarrow}\M_1\stackrel{\iota_1}{\longleftarrow}\M_2
\stackrel{\iota_2}{\longleftarrow}\M_3\stackrel{\iota_3}{\longleftarrow}\cdots.\]
Let $f_\infty:\lim_\iota\F\arr\lim_\iota\F$ be the map induced by
the coverings $f_n$. Then the dynamical system $(\lim_\iota\F,
f_\infty)$ depends, up to topological conjugacy, only on the
combinatorial equivalence class of the topological automaton $\F$.
\end{theorem}

In fact, we prove that the dynamical system $(\lim_\iota\F,
f_\infty)$ is topologically conjugate with the \emph{limit
dynamical system} of the iterated monodromy group of $\F$. Limit
dynamical systems of contracting self-similar groups (contracting
virtual endomorphisms) were defined in~\cite{bgn,nek:book} using
symbolic dynamics (as quotients of the space of infinite sequences
by an equivalence relation defined by a group action).

The inverse limit $\lim_\iota\F$ is an analogue of the Julia set
of an expanding dynamical system. For example, suppose that
$f\in\C(z)$ is a hyperbolic rational function of one complex
variable acting on the Riemann sphere. Then there exists a compact
neighborhood $\M$ of the Julia set of $f$ such that
$f^{-1}(\M)\subset\M$ and $\M$ does not contain the critical
values of $f$. Consider the topological automaton $f,
\iota:\M_1\arr\M$, where $\M_1=f^{-1}(\M)$ and $\iota:\M_1\arr\M$
is the identical embedding. Then the inverse limit
$\lim_\iota\M_n=\bigcap_{n\ge 1}\M_n$ is the Julia set of $f$. The
automaton $(\M, \M_1, f, \iota)$ is contracting with respect to
the restriction onto $\M$ of the Poincar\'e metric on the sphere
minus the post-critical set of $f$.

A partial case of Theorem~\ref{th:main1} (when $\M$ is a
Riemannian manifold and $\iota$ is a diffeomorphism) is the
theorem of M.~Shub on expanding endomorphisms of manifolds,
see~\cite{shub1,shub2}.

Theorem~\ref{th:main1} can be used now to approximate dynamical
systems (acting on their Julia sets) by topological automata. For
instance, if $f:\M_1\arr\M$ is an expanding partial self-covering,
then we can replace $\M, \M_1$, $f$ and the embedding
$\M_1\hookrightarrow\M$ by homotopically equivalent spaces and
maps $f', \iota':\M_1'\arr\M'$, thus getting a topological
automaton $\F$ combinatorially equivalent to the partial
self-covering $f$. If we find a length structure on $\M'$ such
that $\iota':\M_1'\arr\M'$ is contracting with respect to the lift
of the length structure of $\M'$ to $\M_1'$ by $f$, then
Theorem~\ref{th:main1} implies that the dynamical system
$(\lim_{\iota'}\F, f_\infty')$ is topologically conjugate to the
action of $f$ on its Julia set. (Here the Julia set of an
expanding map is defined as the limit set of inverse iterations.)
In particular, the spaces $\M_n'$ approximate the Julia set of
$f$. A known example of this approach are the classical Hubbard
trees of post-critically finite polynomials
(see~\cite{DH:orsayI,DH:orsayII}), which are constructed by
retracting the \emph{Thurston orbispace} of the polynomial onto a
finite tree. Our method has no restrictions on dimension of the
spaces, and can be applied to any expanding dynamical system. For
example, we construct polyhedral models of the Julia sets of
post-critically finite endomorphisms of $\CP^n$ coming from
Teim\"uller theory of post-critically finite polynomials.

A natural question arises now in connection with
Theorem~\ref{th:main1}. How to construct a simple contracting
topological automaton (e.g., consisting of simplicial complexes
$\M, \M_1$ and piecewise affine maps $f, \iota$) with given
iterated monodromy group (with given associated virtual
endomorphism)? Such a construction will provide approximations of
the Julia sets of expanding dynamical systems in a general and
systematic way.

Let $\phi:G_1\arr G$ be a surjective virtual endomorphism of a
finitely generated group $G$. Suppose that $\X$ is a path
connected metric space on which $G$ acts by isometries properly
and co-compactly. Then the identity map on $\X$ induces a covering
$f:\X/G_1\arr\X/G$ of the corresponding orbispaces (if the action
of $G$ on $\X$ is free, then $f$ is a covering of topological
spaces). Suppose that a map $F:\X\arr\X$ is such that
\begin{equation}\label{eq:Fmodel}F(\xi\cdot g)=F(\xi)\cdot\phi(g)\end{equation}
for all $g\in G_1$ and $\xi\in\X$. Then $F$ induces a continuous
map (a morphism of orbispaces) $\iota:\X/G_1\arr\X/G$. We get in
this way a topological automaton $\F=(\X/G_1, \X/G, f, \iota)$. If
$\X$ is simply connected, then $\pi_1(\X/G)=G$ and the virtual
endomorphism associated with the constructed topological automaton
is $\phi$. In general, if $\wt\phi$ is the virtual endomorphism of
$\pi_1(\X/G)$ associated with the automaton $\F$, then $G$ is the
quotient of $\pi_1(\X/G)$ by a normal subgroup invariant under
$\wt\phi$, and $\phi$ is the virtual endomorphism induced by
$\wt\phi$ on the quotient. Iteration of the automaton $\F$
produces the spaces $\M_n=\X/\dom\phi^{\circ n}$. The maps
$f_n:\M_{n+1}\arr\M_n$ are the coverings induced by the inclusions
$\dom\phi^{\circ(n+1)}\le\dom\phi^{\circ n}$; the maps
$\iota_n:\M_{n+1}\arr\M_n$ are induced by the map $F$.

It follows that the question of finding a contracting topological
automaton with given iterated monodromy group is equivalent to the
question of finding a proper co-compact $G$-space $\X$ and a
contracting map $F:\X\arr\X$ satisfying~\eqref{eq:Fmodel}.

The most natural proper co-compact $G$-space is the group $G$
itself with respect to right translations. Choose a left coset
transversal $\{r_1=1, r_2, \ldots, r_d\}$ for the subgroup $G_1$.
Then we can define a map $F:G\arr G$ satisfying~\eqref{eq:Fmodel}
by the formula \[F(g)=\phi(r_ig),\] where $r_i$ is such that
$r_ig\in G_1$. But we need to have a metric space $\X$ such that
$F$ is contracting. A standard approach in geometric group theory
is to consider the Cayley, or Rips complex of $G$
(see~\cite{gro:hyperb}). If $S$ is a finite generating set of $G$,
then denote by $\Gamma(G, S)$ the simplicial complex with the set
of vertices $G$ in which a subset $A\subset G$ is a simplex if and
only if $A\cdot g^{-1}\subset S$ for all $g\in A$. If $S$ is
invariant with respect to the map $F$, then we get a simplicial
map $F:\Gamma(G, S)\arr\Gamma(G, S)$ satisfying~\eqref{eq:Fmodel}.

It is proved in Subsection~\ref{ss:contractingautomaton}
(Theorem~\ref{th:contractingautomaton}) that this natural
construction works.

\begin{theorem}
\label{th:main2} If $\phi:G_1\arr G$ is a contracting virtual
endomorphism (e.g., the virtual endomorphism associated with a
contracting automaton), then there exist positive integers $m$ and
$n$ such that the map $F^{\circ n}:\Gamma(G, S^m)\arr\Gamma(G,
S^m)$ is homotopic through maps satisfying~\eqref{eq:Fmodel} to a
contracting map.
\end{theorem}

In this way we get for every contracting topological automaton
$\mathcal{F}=(\M, \M_1, f, \iota)$ (e.g., for every expanding
partial self-covering) a contracting simplicial topological
automaton combinatorially equivalent to some iteration of
$\mathcal{F}$. The Julia set of $\mathcal{F}$ will be homeomorphic
to the inverse limit of the simplicial complexes $\Gamma(G,
S^m)/\dom\phi^{\circ nk}$ as $k\to\infty$.

Note that every finite-dimensional compact metric space is an
inverse limit of simplicial complexes, by an old theorem of
P.~Alexandroff~\cite{alexandroff:simplicialapprox}.

Theorem~\ref{th:contractingautomaton} proved in our paper is more
explicit and ``cleaner'' than Theorem~\ref{th:main2}. In
particular, we use a smaller simplicial complex than $\Gamma(G,
S^m)$ (the map $F$ is not surjective on $\Gamma(G, S^m)$, so we
can pass to the intersection of domains of its iteration).

Due to combinatorial nature of the simplicial complex $\Gamma(G,
S^m)$, the complexes $\Gamma(G, S^m)/\dom\phi^{\circ n}$ are
constructed using simple recursive cut-and-paste rules, described
in Proposition~\ref{pr:cutandpaste}.

The structure of the paper is as follows. The second section is an
overview of the techniques of self-similar groups, virtual
endomorphisms, and their limit spaces. All proofs can be found in
the monograph~\cite{nek:book}.

In Section~\ref{s:topautomata} we define topological automata and
describe some examples (Moore diagrams of finite automata, wreath
recursions, partial self-coverings, post-critically finite
rational functions, post-critically finite correspondences,
bi-reversible automata, commensurizers of tree lattices, Thurston
maps, and subdivision rules).

In Section~\ref{s:img} we show how topological automata are
iterated; define the inverse limit $\lim_\iota\F$, and two other
inverse limits $\lim_f\F$ and $\lim_{f, \iota}\F$; define iterated
monodromy groups of topological automata; and show how they are
computed as self-similar groups. At the end of the section we
define the notion of combinatorial equivalence of topological
automata.

Section~\ref{s:contracting} studies contracting topological
automata. We pass to a more convenient setting of group actions on
topological spaces. If $(\M, \M_1, f, \iota)$ is a topological
automaton, then passing to the universal covering $\X$ of $\M$ we
get an action of $\pi_1(\M)$ on $\X$, a subgroup
$G_1\cong\pi_1(\M_1)$ of $\pi_1(\M)$, and a map $F:\X\arr\X$,
which is a lift of the map $\iota$ to the universal covering. The
map $F$ satisfies the condition~\eqref{eq:Fmodel} for the virtual
endomorphism $\phi$ associated with the topological automaton. We
formalize such structures, and pass from the study of topological
automata to the study of group actions and equivariant maps. We
prove then results equivalent to Theorem~\ref{th:main1}: one is
formulated in terms of group actions
(Theorem~\ref{th:approximationlimg}), and the other in terms of
topological automata (Theorem~\ref{th:projlim}). We also show that
a topological automaton homotopy equivalent to a contracting
topological automaton can be made contracting, if we pass to its
iteration (Corollary~\ref{cor:hequiv}). This result can be used to
construct contracting topological automata approximating an
expanding dynamical system by passing to homotopy equivalent
spaces and maps.

In Section~\ref{s:simplicialapprox} we show how to construct a
contracting piecewise affine topological automaton starting from
any contracting iterated monodromy group (i.e., starting from any
contracting virtual endomorphism of a group). Our construction
essentially coincides with the one described in
Theorem~\ref{th:main2} above. The only difference is that we pass
to the smaller complex $\bigcap_{k\ge 1}F^k(\Gamma(G, S^m))$,
which will not depend now on the choice of $S$ and $m$ (if $m$ is
big enough). We also describe recurrent cut-and-paste rules for
constructing the simplicial complexes approximating the Julia set
(Proposition~\ref{pr:cutandpaste}).

The last section presents some examples of application of the
developed technique. In particular, we show how the Hubbard trees
fit into our theory, and describe polyhedral models of
post-critically finite rational endomorphisms of complex
projective spaces coming from Teichm\"uller theory of hyperbolic
polynomials.

\medskip
\noindent\textbf{Acknowledgements.} The author is very grateful
for fruitful discussions with Laurent Bartholdi, Andr\'e
Haefliger, Sarah Koch, and John Smillie on the topics of these
notes.

\section{Self-similar groups and their limit spaces}
We give in this section a short overview of the main definitions
and constructions of the theory of self-similar groups. For more
details and proofs, see~\cite{nek:book,nek:filling}.

\subsection{Main definitions}
For a finite set $\alb$, we denote by $\xs=\bigsqcup_{n\ge
0}\alb^n$ the set of finite words over $\alb$, i.e., the free
monoid generated by $\alb$.

\begin{defi}
A \emph{faithful self-similar action} $(G, \alb)$ is a faithful
action of a group $G$ on the set $\xs$ such that for every $g\in
G$ and $x\in\alb$ there exist $h\in G$ such that
\[g(xw)=g(x)h(w)\]
for all $w\in\xs$.
\end{defi}

Every self-similar action preserves the levels $\alb^n$ of $\xs$.
It follows from the definition that for every word $v\in\xs$ and
every $g\in G$ there exists $h\in G$ such that
\[g(vw)=g(v)h(w).\]
The element $h$ is unique, by faithfulness of the action. We
denote $h=g|_v$ and call $h$ the \emph{section} (or
\emph{restriction}) of $g$ at $v$. We have the following
properties of sections
\begin{equation}\label{eq:sections}
g|_{v_1v_2}=g|_{v_1}|_{v_2},\quad(g_1g_2)|_v=g_1|_{g_2(v)}g_2|_v.
\end{equation}

Self-similar actions are usually described by the \emph{associated
wreath recursion}, which is the homomorphism
$\varphi:G\arr\symm\wr G=\symm\ltimes G^\alb$ given by
\[\varphi(g)=\pi\cdot(g|_x)_{x\in\alb},\]
where $\pi$ is the permutation of $\alb=\alb^1\subset\xs$ defined
by $g$.

For example, the transformation of $\{0, 1\}^*$ defined by the
recursive rules
\[a(0w)=1w,\quad a(1w)=0a(w)\]
is defined in terms of the associated wreath recursion as
\[\varphi(a)=\sigma(1, a),\]
where $\sigma\in\symm[\{0, 1\}]$ is the transposition and $1$ on
the right hand side of the equality is the trivial transformation.
We will usually omit  $\varphi$ and write the last equality as
$a=\sigma(1, a)$.

The elements of the wreath product $\symm\ltimes G^\alb$ are
multiplied according to the rule
\[\pi_1(g_x)_{x\in\alb}\cdot\pi_2(h_x)_{x\in\alb}=\pi_1\pi_2(g_{\pi_2(x)}h_x)_{x\in\alb}.\]
Note that we use left action in this formula.

The wreath recursion uniquely determines the associated
self-similar action. The following proposition is proved
in~\cite[Proposition~2.3.4]{nek:book} (see
also~\cite[Proposition~2.12]{nek:filling}).

\begin{proposition}
Let $(G, \alb)$ be a self-similar action and let
$\phi:G\arr\symm\wr G$ be the associated wreath recursion. For
every element $h\in\symm\wr G$ the self-similar action associated
with the wreath recursion $g\mapsto h^{-1}\phi(g) h$ is conjugate
to the self-similar action $(G, \alb)$.
\end{proposition}

\begin{defi}
Two self-similar actions of a group $G$ on $\xs$ are said to be
\emph{equivalent} if the associated wreath recursions can be
obtained from each other by taking composition with an inner
automorphism of the group $\symm\wr G$.
\end{defi}

An approach equivalent to wreath recursions, but in some sense
more ``coordinate-free'', uses the notion of a \emph{permutational
bimodule}, which is defined as follows.

\begin{defi}
Let $G$ be a group. A \emph{permutational $G$-bimodule} is a set
$\bim$ together with commuting left and right actions of $G$ on
it, i.e., with two maps $G\times\bim\arr\bim:(g, x)\mapsto g\cdot
x$ and $\bim\times G:(x, g)\mapsto x\cdot g$ satisfying the
following conditions
\begin{enumerate}
\item $1\cdot x=x$, $x\cdot 1=x$ for all $x\in\bim$;
\item $g_1\cdot(g_2\cdot x)=(g_1g_2)\cdot x$ and $(x\cdot
g_1)\cdot g_2=x\cdot(g_1g_2)$ for all $g_1, g_2\in G$ and
$x\in\bim$;
\item $g_1\cdot(x\cdot g_2)=(g_1\cdot x)\cdot g_2$ for all $g_1,
g_2\in G$.
\end{enumerate}
\end{defi}

Let $(G, \alb)$ be a self-similar action. If we identify the
letters of $\alb$ with the transformations
\[x:v\mapsto xv\]
of $\xs$, then the condition
\[g(xw)=yh(w)\quad\forall w\in\xs\]
is written as the equality
\[g\cdot x=y\cdot h\]
of compositions of transformations of $\xs$. It follows that the
set $\alb\cdot G$ of transformations $x\cdot g$ for $x\in\alb$ and
$g\in G$ is a $G$-bimodule with respect to pre- and
post-compositions with the elements of $G$. The obtained bimodule
is called the \emph{associated bimodule} of the self-similar
action (or the \emph{self-similarity bimodule}).
The left and right actions of $G$ on the set $\alb\cdot G$
are given then by the rules
\[h\cdot (x\cdot g)=h(x)\cdot (h|_xg),\quad (x\cdot g)\cdot
h=x\cdot (gh).\] The right action of $G$ on $\alb\cdot G$ is free
(i.e., $x\cdot g=x$ implies $g=1$) and has $|\alb|$ orbits. We
generalize these conditions in the following definition.

\begin{defi}
\label{def:covbim} A \emph{($d$-fold) covering $G$-bimodule} is a
permutational $G$-bimodule $\bim$ such that the right action of
$G$ on $\bim$ is free and has $d$ orbits.

A transversal $\alb\subset\bim$ of the right orbits, i.e., a set
intersecting every orbit of the right action exactly once, is
called a \emph{basis} of the covering bimodule $\bim$.
\end{defi}

Let $\bim$ be a $d$-fold covering $G$-bimodule. Choose a basis
$\alb$. Then every element of $\bim$ can be uniquely written in
the form $x\cdot g$ for $x\in\alb$ and $g\in G$. Consequently, for
every $g\in G$ and $x\in\alb$ there exist unique $h\in G$ and
$y\in\alb$ such that $g\cdot x=y\cdot h$ in $\bim$. The
\emph{associated self-similar action} $(G, \alb, \bim)$ of $G$ on
$\xs$ is given then by the recurrent rule
\[g(xw)=yh(w)\Longleftarrow g\cdot x=y\cdot
h.\]

The action $(G, \alb, \bim)$ does not depend, up to equivalence of
the actions (hence up to conjugacy), on the choice of the basis
$\alb$.

The action associated to a covering $G$-bimodule $\bim$ is not
faithful in general. The \emph{faithful quotient} of $G$ is the
quotient of $G$ by the kernel of the associated action. The action
of the faithful quotient on $\xs$ is self-similar and the
associated bimodule is called the \emph{faithful quotient} of the
bimodule $\bim$.

If $\bim_1$ and $\bim_2$ are permutational $G$-bimodules, then
their tensor product $\bim_1\otimes\bim_2$ is the quotient of
$\bim_1\times\bim_2$ by the identifications $x_1\cdot g\otimes
x_2=x_1\otimes g\cdot x_2$. It is a $G$-bimodule with respect to
the actions $g_1\cdot (x_1\otimes x_2)\cdot g_2=(g_1\cdot
x_1)\otimes (x_2\cdot g_2)$. If $\bim_1$ and $\bim_2$ are covering
bimodules, then $\bim_1\otimes\bim_2$ is also a covering bimodule.

If $\bim$ is a covering bimodule and $\alb$ is its basis, then the
set $\alb^n$ of words $x_1x_2\ldots x_n=x_1\otimes
x_2\otimes\cdots\otimes x_n$, for $x_i\in\alb$, is a basis of the
bimodule $\bim^{\otimes n}$. For every $v\in\alb^n$ and $g\in G$
there exists then a unique pair $u\in\alb^n$ and $h\in G$ such
that $g\cdot v=u\cdot h$ in $\bim^{\otimes n}$. The action
$g:v\mapsto u$ coincides then with the associated self-similar
action $(G, \alb, \bim)$.

\begin{defi}
A \emph{virtual endomorphism} of a group $G$ is a homomorphism of
groups $\phi:G_1\arr G$, where $G_1<G$ is a subgroup of finite
index (called the \emph{domain} of $\phi$).

Two virtual endomorphisms $\phi_1, \phi_2$ of $G$ are
\emph{conjugate} if there exist $g, h\in G$ such that
$h^{-1}\cdot\dom\phi_1\cdot h=\dom\phi_2$ and
$\phi_1(x)=g^{-1}\phi_2(h^{-1}xh)g$ for all $x\in\dom\phi_1$.
\end{defi}

If $\bim$ is a covering $G$-bimodule then, for $x\in\bim$, the
\emph{associated virtual endomorphism} $\phi_x$ is given by the
rule
\[g\cdot x=x\cdot\phi_x(g),\]
and is defined on the subgroup of the elements $g\in G$ such that
$x$ and $g\cdot x$ belong to one right $G$-orbit. If the left
action of $G$ on the set of right orbits is transitive, then the
bimodule $\bim$, and the associated self-similar action are
uniquely determined (up to isomorphism of the bimodules and up to
equivalence of self-similar actions) by the associated virtual
endomorphism (see~\cite[Proposition~2.5.8]{nek:book}).

\subsection{Contracting groups and their limit spaces}
\label{ss:contractinglimsp}
\begin{defi}
A self-similar group $(G, \alb)$ is said to be \emph{contracting}
if there exists a finite set $\nuke\subset G$ such that for every
$g\in G$ there exists $n_0\in\N$ such that $g|_v\in\nuke$ for all
words $v\in\xs$ of length at least $n_0$.
\end{defi}

If the group is contracting, then the smallest set $\nuke$
satisfying the conditions of the definition is called the
\emph{nucleus} of the action.

If a self-similar group is contracting, then every equivalent
action is also contracting (though the nucleus may be different).
Consequently, the property of being contracting depends only on
the associated bimodule, and does not depend on the choice of the
basis.

Denote by $\xmo$ the set of left-infinite sequences $\ldots
x_2x_1$ over the alphabet $\alb$ with the direct product topology
(where $\alb$ is discrete).

\begin{defi}
Let $(G, \alb)$ be a contracting group. We say that two sequences
$\ldots x_2x_1, \ldots y_2y_2\in\xmo$ are \emph{asymptotically
equivalent} with respect to the action $(G, \alb)$ if there exists
a finite set $N\subset G$ and a sequence $g_k\in N$ such that
\[g_k(x_k\ldots x_2x_1)=y_k\ldots y_2y_1\]
for all $k$. The quotient of the space $\xmo$ by the asymptotic
equivalence relation is called the \emph{limit space} of the
action and is denoted $\lims$.
\end{defi}

\begin{proposition}
The limit space of a contracting self-similar group is a finite
dimensional compact metrizable space. The shift $\ldots
x_2x_1\mapsto\ldots x_3x_2$ on $\xmo$ induces a continuous map
$\si:\lims\arr\lims$.
\end{proposition}

The dynamical system $(\lims, \si)$ is called the \emph{limit
dynamical system} of the contracting group $(G, \alb)$.

A natural structure of an \emph{orbispace} on $\lims$ is
introduced using the following ``covering space'' of $\lims$.

\begin{defi}
Let $(G, \alb)$ be a contracting self-similar group. Let
$\xmo\times G$ be the direct product of the topological space
$\xmo$ with the discrete group $G$. We say that $\ldots
x_2x_1\cdot g$ and $\ldots y_2y_1\cdot h\in\xmo\times G$ are
\emph{asymptotically equivalent} if there exists a finite set
$N\subset G$ and a sequences $g_k\in G$, such that
\[g_k(x_k\ldots x_2x_1)=(y_k\ldots y_2y_1),\quad g_k|_{x_k\ldots
x_2x_1}g=h\] for all $k$. The quotient of $\xmo\times G$ by the
asymptotic equivalence relation is called the \emph{limit
$G$-space} and is denoted $\limg$.
\end{defi}

It is easy to see that the natural right action of $G$ on
$\xmo\times G$ induces an action of $G$ on $\limg$. The space of
orbits $\limg/G$ of this action is homeomorphic to $\lims$. The
corresponding orbispace is the \emph{limit orbispace} of the
contracting group $(G, \alb)$. For theory of orbispaces,
see~\cite[Chapter~III.$\mathcal{G}$]{bridhaefl}
and~\cite[Chapter~4]{nek:book}.

For $x\in\alb$ and a point $\xi\in\limg$ represented by a sequence
$\ldots x_2x_1\cdot g$ we denote by $\xi\otimes x$ the point
represented by $\ldots x_2x_1g(x)\cdot g|_x$. The map $\xi\mapsto
\xi\otimes x$ is continuous.

\begin{defi}
Let $(G, \alb)$ be a contracting group. The \emph{tile} $\til$ is
the image of the set $\xmo\cdot\{1\}\subset\xmo\times G$ in the
limit $G$-space $\limg$.

For $v\in\alb^n$ and $g\in G$, the corresponding \emph{tile of
$n$th level} is the set $\til\otimes v\cdot g$, i.e., the image in
$\limg$ of the set of sequences ending by $v\cdot g$.
\end{defi}

\begin{proposition}
\label{pr:tilesintersection} Two tiles $\til\otimes v_1\cdot g_1$
and $\til\otimes v_2\cdot g_2$ of the $n$th level intersect if and
only if there exists an element $h$ of the nucleus such that
$h\cdot v_1\cdot g_1=v_2\cdot g_2$.
\end{proposition}

We write $h\cdot v=u\cdot g$, for $v, u\in\alb^n$ and $h, g\in G$,
if $h(v)=u$ and $g=h|_v$ (which agrees with the definition of the
bimodule $\bim^{\otimes n}$). In particular, the equality $h\cdot
v_1\cdot g_1=v_2\cdot g_2$ means that $v_2=h(v_1)$ and
$h|_{v_1}g_1=g_2$.

\section{Topological automata}
\label{s:topautomata}
\subsection{Definition}
\begin{defi}
A \emph{topological automaton} is a quadruple $\F=(\M, \M_1, f,
\iota)$, where $\M$ and $\M_1$ are topological spaces (or
orbispaces), $f:\M_1\arr\M$ is a finite covering map and
$\iota:\M_1\arr\M$ is a continuous map (a morphism of orbispaces).
\end{defi}

This definition coincides (in the regular, i.e., non-orbispace
case) with the notion of a \emph{topological correspondence} or
\emph{topological graph} studied by T.~Katsura
in~\cite{katsura:generalizing,katsura:examples,katsura:ideals,katsura:pureinfiniteness},
which might be a better terminology. We use a different term in
order to show a strong connection to the theory of self-similar
groups and groups generating by automata. We will consider
topological automata up to different weak equivalence relations,
so that they will be combinatorial rather than rigidly topological
objects.

\begin{defi}
\label{def:thurstequiv} Topological automata $\F=(\M, \M_1, f,
\iota)$ and $\F'=(\M', \M_1', f', \iota')$ are said to be
\emph{homotopy
  equivalent} if there exist homotopy equivalences
$\phi_1:\M_1'\arr\M_1$, $\phi:\M'\arr\M$ and maps
$\iota_1':\M_1'\arr\M'$, $\iota_1:\M_1\arr\M$, homotopic to
$\iota'$ and $\iota$, respectively, such that the diagrams
\[\begin{array}{ccc}\M_1' & \stackrel{\phi_1}{\arr} & \M_1\\
\mapdown{f'} & & \mapdown{f}\\
\M' & \stackrel{\phi}{\arr} & \M\end{array}\qquad
\begin{array}{ccc}\M_1' & \stackrel{\phi_1}{\arr} & \M_1\\
\mapdown{\iota'_1} & & \mapdown{\iota_1}\\
\M' & \stackrel{\phi}{\arr} & \M\end{array}\] are commutative.
\end{defi}

We consider topological automata up to homotopy equivalence. We
will introduce an even weaker equivalence relation between
topological automata later.

Here topological automata are topological analogs of transducers.
They should not be confused with analogs of \emph{acceptors} (see,
for instance~\cite{brauer:french,jeandel:topautom}).

\subsection{Examples of topological automata}

\subsubsection{Automata and Moore diagrams}
\label{sss:moorediagrams} Let us recall the definition of
\emph{automata} (also known as \emph{transducers}). For more on theory
of transducers and groups generated by automata,
see~\cite{eil,grineksu_en}.

\begin{defi}
An \emph{automaton} over the alphabet $\alb$ is a triple $(Q, \pi,
\tau)$, where $Q$ is a set (of \emph{internal states}), and $\pi$
and $\tau$ are maps
\[\pi:Q\times\alb\arr\alb,\quad\tau:Q\times\alb\arr Q,\]
called the \emph{output} and \emph{transition} functions,
respectively. The automaton is called \emph{invertible} if for
every $q_0\in Q$ the map $x\mapsto\pi(q_0, x)$ is a permutation.
The automaton is \emph{finite} if the set $Q$ is finite.
\end{defi}

We interpret the automaton $(Q, \pi, \tau)$ as a machine, which
being in a state $q\in Q$ and reading on input a letter $x$ prints
the letter $\pi(q, x)$ on the output and changes its state to
$\tau(q, x)$.

Every invertible automaton can be related to a topological
automaton by the notion of a \emph{Moore diagram} (also called a
\emph{state diagram}).

\begin{figure}
\centering\includegraphics{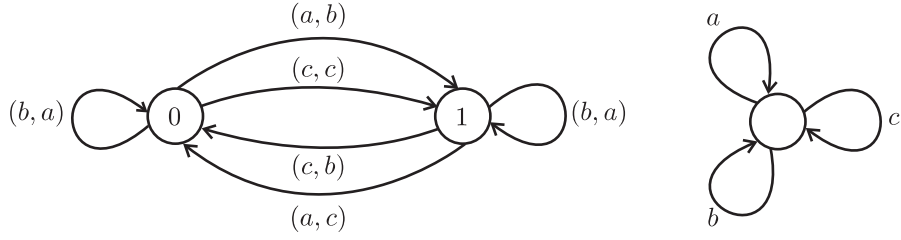}\\
\caption{Dual Moore diagram}\label{fig:dualMoore}
\end{figure}

Moore diagrams are classical representations of automata. We will
use here the \emph{dual Moore diagrams} (i.e., the usual Moore
diagrams of the dual automata). It is an oriented graph with the
set of vertices $\alb$ in which for every $x\in\alb$ and $q\in Q$
we have an arrow starting in $x$, ending in $\pi(q, x)$ and
labeled by $(q, \tau(q, x))$. See an example of a dual Moore
diagram of an automaton on Figure~\ref{fig:dualMoore}.

Dual Moore diagrams are naturally interpreted as topological
automata. We take $\M_1$ to be the dual Moore diagram of the
automaton $(Q, \pi, \tau)$ as a topological graph (i.e., as a
cellular complex). The space $\M$ is a graph with one vertex and
$|Q|$ loops labeled by the elements of $Q$. If an arrow of $\M_1$
is labeled by $(q_1, q_2)$, then it is mapped by $f:\M_1\arr\M$ to
the loop of $\M$ labeled by $q_1$ and by $\iota:\M_1\arr\M$ to the
loop labeled by $q_2$. We get in this way a topological automaton,
which is well defined up to a homotopy equivalence. The condition
of invertibility of the automaton $(Q, \tau, \pi)$ is equivalent
to the condition that $f$ is a covering map.

\subsubsection{Wreath recursions}
More generally, let $(G, \alb)$ be a finitely generated
self-similar group. Choose a generating set $S$ of $G$. Consider
the \emph{dual Moore diagram} of $(G, \alb)$ with respect to $S$.
It is a directed graph with the set of vertices $\alb$ in which
for every $g\in S$ and $x\in\alb$ we have an arrow starting in
$x$, ending in $g(x)$, and labeled by $(g, g|_x)$. This graph
describes the wreath recursion of $(G, \alb)$: the arrows labeled
by $(g, \cdot)$ describe the action of the generator $g$ on
$\alb$, and the second coordinates of the labels show the
corresponding sections $g|_x$.

The dual Moore diagram of a self-similar group is a topological
automaton. The graph $\M$, as in the previous example, has one
vertex and oriented loops labeled by the elements of $S$. Let
$\M_1$ be the dual Moore diagram of $(G, \alb)$. The first
coordinates of the labels show the values of the covering
$f:\M_1\arr\M$; the second coordinates show the values of the map
$\iota:\M_1\arr\M$: the arrow labeled by $(g, h)$ is mapped by
$\iota$ to the path in $\M$ such that the product of the labels
along the path (taking into account the orientation) is equal to
$h$. Note that the obtained topological automaton is not uniquely
defined even up to a homotopy equivalence, since elements $h$ of
$G$ may be represented in different ways as products of the
elements of $S$.

\subsubsection{Partial self-coverings} If $\iota$ is an embedding,
then the topological automaton $(\M, \M_1, f, \iota)$ is a
\emph{partial self-covering} of $\M$. Partial self-coverings are
studied in~\cite{nek:book}. See also~\cite{nek:filling}, where the
category of partial self-coverings is defined. Theory of
topological automata and their iterated monodromy groups is not
much different from the theory of partial self-coverings. The main
reason to introduce the general notion (except for the pure sake
of generality) is that topological automata are less rigid
objects, and are easier to construct, and hence to use them as
models of more complicated partial self-coverings and their Julia
sets.

\subsubsection{Post-critically finite rational functions}
A rational function $f:\CS\arr\CS$ is said to be
\emph{post-critically finite} if the orbit of every critical point
of $f$ under iterations of $f$ is finite. Denote by $P$ the
\emph{post-critical set} of $f$, i.e., the union of the orbits of
critical values of $f$. Then $f:\CS\setminus
f^{-1}(P)\arr\CS\setminus P$ is a partial self-covering (since
$\CS\setminus f^{-1}(P)\subseteq\CS\setminus P$). Hence
post-critically finite rational functions are examples of
topological automata.

\subsubsection{Post-critically finite correspondences}\label{sss:correspond} Let
$R\subset\CS\times\CS$ be a correspondence, i.e., an algebraic
curve in $\CS\times\CS$. Denote by $p_1$ and $p_2$ projections
of $R$ onto the first and the second coordinates of the
correspondence. We assume that $p_1$ and $p_2$ are branched
coverings and interpret $R$ as a multivalent function
\[p_2(z)\mapsto p_1(z).\]

Suppose that $R$ is \emph{post-critically finite}, i.e., there
exists a finite set $P\subset R$ such that $p_1:R\setminus
P\arr\CS\setminus p_1(P)$ is a covering and $p_1(P)\subseteq
p_2(P)$.

We have $p_2(R\setminus P)\subseteq p_1(R\setminus P)$, hence the
quadruple $(p_1(R\setminus P), R\setminus P, p_1, p_2)$ is a
topological automaton.

As a simple example, consider the correspondence
\[z^q\mapsto z^p\]
for natural numbers $p$ and $q$, i.e., the multivalent function
$z^{p/q}$. Its post-critical set is $\{0, \infty\}$.

Another famous example is the correspondence associated with the
\emph{arithmetic-geometric mean}, studied by Gauss. An extensive
account on the history of arithmetic-geometric mean is given
in~\cite{cox:agm}.

Lagrange in 1784 and independently Gauss in~1790 have shown that
if $a_0$ and $b_0$ are positive real numbers, then the sequences
\[a_n=\frac 12(a_{n-1}+b_{n-1}),\quad b_n=\sqrt{a_{n-1}b_{n-1}}\]
converge to a common value $M(a_0, b_0)$, called the
arithmetic-geometric mean.

In the complex case one has to choose one of two signs of the
square root. We get the correspondence
\[[z_1:z_2]\mapsto\left[(z_1+z_2)/2:\sqrt{z_1z_2}\right],\]
on the projective line $\CS$. It is written in the affine
coordinates as
\[w\mapsto\frac{1+w}{2\sqrt{w}}.\]

In our terms, the correspondence is given by the pair of maps
\[f(w)=\frac{(1+w)^2}{4w},\quad\iota(w)=w^2,\]
so that it is the curve $\left\{\left(\frac{(1+w)^2}{4w},
w^2\right)\;:\;w\in\CS\right\}$. Denote by $P$ the set of the
points $(\infty, 0), (1, 1), (0, 1)$ and $(\infty, \infty)$, which
are the points of $R$ parametrized by $w=0, 1, -1$ and $\infty$,
respectively. We have
\[f(P)=\{\infty, 1, 0\}=\iota(P),\] and the maps $f,
\iota:R\setminus P\arr \CS\setminus\pi_1(P)$ are coverings.

See~\cite{bullett:agm}, where the arithmetic-geometric mean is
studied as an example of a post-critically finite correspondence.
For more on dynamics of correspondences see the
papers~\cite{bullett:quadratic,bullet:crfinitemodular,bullett:gallery}.

\subsubsection{Bi-reversible automata}
In the last example both maps $f$ and $\iota$ were coverings. This
situation for the dual Moore diagrams of automata has a special name.

\begin{defi}
Let $(\M, \M_1, f, \iota)$ be the dual Moore diagram of a finite
invertible automaton. The automaton is said to be
\emph{bi-reversible} if the map $\iota$ is a covering.
\end{defi}

An example of a bi-reversible automaton (of its dual Moore
diagram) is shown on Figure~\ref{fig:dualMoore}. It corresponds to one of two
automata, which appeared in the paper~\cite{al:free_en}. It was
proved in~\cite{vorobets:alfree} that the self-similar group
generated by this automaton is free. For more on bi-reversible
automata see~\cite{mns_en,glasnermozes,vorobets:alfree} and
Section~1.10 of~\cite{nek:book}.

\subsubsection{Commensurizers of tree lattices}
If $G$ is a group and $H<G$ is a subgroup, then the
\emph{commensurizer} of $H$ in $G$ is the group of the elements
$g\in G$ such that $H\cap g^{-1}Hg$ has finite index in $H$ and
$g^{-1}Hg$.

As a generalization of bi-reversible automata, consider the
topological automata $(\M, \M_1, f, \iota)$, where $\M$ is a
bouquet of $k$ circles, and $f:\M_1\arr\M$ and $\iota:\M_1\arr\M$
are coverings. It is shown in Proposition~2.2 of~\cite{lmz} that
the topological automata $(\M, \M_1, f, \iota)$ of this form
describe the elements of the commensurizer of the co-compact
lattice $\pi_1(\M)$ in the automorphism group of the universal
covering $T$ of $\M$.

More precisely, if $g$ is an element of the commensurizer of
$\pi_1(\M, t)$ in the automorphism group of $T$, then there exists
a finite index subgroup $H<\pi_1(\M)$ such that $f(e)=\iota(g(e))$
for every edge $e$ of $T$. It follows that $g$ is uniquely
determined by the image $g(t_0)$ of a vertex $t_0$ of $T$ and by
the automaton $(\M, \M_1, f, \iota)$ (which is called a
\emph{periodic recoloring} in~\cite{lmz}).

For more on lattices in the automorphism groups of trees,
see~\cite{lub:treelatices,glasnermozes}.

\subsubsection{Thurston maps}
\label{sss:thurstonmaps} A \emph{Thurston map} is a
post-critically finite orientation preserving branched
self-covering $f:S^2\arr S^2$ of the sphere. It can be interpreted
as a topological automaton in the same way as in the case of
post-critically finite rational functions.

Thurston's theorem (see~\cite{DH:Thurston}) gives a criterion when
a Thurston map is equivalent to a post-critically finite rational
function. Definition of homotopy equivalence
(Definition~\ref{def:thurstequiv}) is a generalization of the
equivalence relation introduced in Thurston's theorem.

In many cases it is more convenient not to remove all
post-critical points from the sphere $S^2$, but rather to
introduce an orbifold structure on $S^2$ minus some post-critical
points. The corresponding orbifold construction, also due to
Thurston, is defined as follows.

Let $C_f$ be the set of critical points of a Thurston map
$f:S^2\arr S^2$ and let $P_f=\bigcup_{n\ge 1}f^n(C_f)$ be the
post-critical set. Let $P'\subset P_f$ be the union of all cycles
of $f$ intersecting $C_f$ (they are \emph{superattracting} if $f$
is a rational function).

The underlying space of the orbifold $\M$ will be $S^2\setminus
P'$. The points $P_f\setminus P'$ will be its singular points.

Denote by $\nu(x)$ for $x\in S^2\setminus P'$ the least common
multiple of the local degrees of $f^m$ at $z$, for all $z$ such
that $f^m(z)=x$. The number $\nu(x)$ is finite for all $x\in
S^2\setminus P'$ and greater than 1 if and only if $x\in P_f$.

Then for any $x\in S^2$ the number $\nu(f(x))$ is divisible by
$\deg_x(f)\cdot \nu(x)$, where $\deg_x(f)$ denotes the local
degree of $f$ at $x$.

Let $\M$ be the orbispace with the underlying space $S^2\setminus
P'$ for which a point $x\in\M$ is uniformized in the atlas of the
orbispace by the cyclic group of order $\nu(x)$ acting by
rotations of a disc.

Similarly, let $\M_1$ be the orbifold defined by the weights
$\nu_0(x)=\frac{\nu(f(x))}{\deg_x(f)}$ instead of $\nu(x)$. The
set of singular points of $\M_1$ is contained in $f^{-1}(P_f)$.
The underlying space of $\M_1$ is $S^2\setminus
f^{-1}\left(P'\right)$.

We have $\nu(z)|\nu_0(z)$, hence the orbispace $\M_{\nu_0}$ is an
open sub-orbispace of $\M_\nu$, where the embedding is identical
on the underlying spaces (see a definition of embedding of
orbifolds in~\cite[Section~4.3]{nek:book}).

On the other hand, the map $f:\M_{\nu_0}\arr\M_{\nu}$ is a
covering of the orbispaces, since
$\deg_x(f)=\frac{\nu_0(x)}{\nu(f(x))}$.

In this way we get an orbifold topological automaton $(\M, \M_1,
f, \iota)$, where $\iota$ is the identical embedding of the
orbispaces.

\subsubsection{Subdivision rules}
Finite subdivision rules are convenient combinatorial descriptions
of Thurston maps, see~\cite{cfp:subdivision,cfp:subdfromrat}.

See~\cite{cfp:subdivision} for a precise definition of subdivision
rules. In our terminology, subdivision rules correspond to
topological automata $(\M, \M_1, f, \iota)$ such that $\M$ and
$\M_1$ are two-dimensional CW complexes (or complexes of groups),
$f:\M_1\arr\M$ is a cellular covering map and $\iota:\M_1\arr\M$
is a homeomorphism such that $\iota(\M_1)$ is a subdivision of
$\M$.

The cells of $\M$ are called \emph{types}. Description of the
covering $f$ amounts to prescribing types (i.e., images under $f$)
to the cells of $\M_1$. The subdivision rule specifies then how
the cells of $\M$ are subdivided into the images of the cells of
$\M_1$ under $\iota$, i.e., specifies the subdivision and labels
the cells according to their types. One also has to label the
edges and vertices appropriately, so that one gets uniquely
defined maps $f$ and $\iota$.

\section{Iterated monodromy groups}
\label{s:img}
\subsection{Iteration of topological automata}
Every topological automaton $\F=(\M, \M_1, f, \iota)$ can be
\emph{iterated} in the following way. Denote $\M_0=\M$, $f_0=f$,
and $\iota_0=\iota$. Define inductively the covering
$f_n:\M_{n+1}\arr\M_n$ as the pullback of the covering
$f_{n-1}:\M_n\arr\M_{n-1}$ by the map
$\iota_{n-1}:\M_n\arr\M_{n-1}$, and the map
$\iota_n:\M_{n+1}\arr\M_n$ as the morphism closing the pullback
diagram \[\begin{array}{ccc}\M_{n+1} & \stackrel{\iota_n}{\arr} &
\M_n\\
\mapdown{f_n} & & \mapdown{f_{n-1}}\\
\M_n & \stackrel{\iota_{n-1}}{\arr} & \M_{n-1}.\end{array}\]

Then the $n$th iteration $\F^n$ of the topological automaton $\F$
is the pair of maps \[f^n=f_0\circ f_1\circ\cdots\circ
f_{n-1},\quad\iota^n=\iota_0\circ\iota_1\circ\cdots\circ
\iota_{n-1}:\M_n\arr\M.\]

In the case when $\M$ and $\M_1$ are regular (i.e., are usual
topological spaces), the pullback $\M_2$ can be defined as the
subspace
\[\{(x, y)\in\M_1^2\;:\;f(y)=\iota(x)\},\]
so that the map $\iota_1:\M_2\arr\M_1$ and the covering
$f_1:\M_2\arr\M_1$ are given by $\iota_1(x, y)=y$ and $f_1(x,
y)=x$.

We get hence by induction the following description of the
iteration.

\begin{proposition}
\label{pr:iteratreg} Let $\F=(\M, \M_1, f, \iota)$ be a
topological automaton such that $\M$ (and hence $\M_1$) are
regular. Then the space $\M_n$ is homeomorphic to the subspace
\[\{(x_1, x_2, \ldots, x_n)\in\M_1^n\;:\;f(x_{i+1})=\iota(x_i),
i=1,\ldots, n-1\},\] and the maps $f_n:\M_{n+1}\arr\M_n$ and
$\iota_n:\M_{n+1}\arr\M_n$ are given by
\[f_n(x_1, x_2, \ldots, x_{n+1})=(x_1, x_2, \ldots, x_n),\]
and
\[\iota_n(x_1, x_2, \ldots, x_{n+1})=(x_2, x_3, \ldots, x_{n+1}).\]
In particular, the topological automaton $\F^n$ is defined by the
maps
\[f^n(x_1, x_2, \ldots, x_n)=f(x_1),\]
and
\[\iota^n(x_1, x_2, \ldots, x_n)=\iota(x_n).\]
\end{proposition}

\begin{examp}
If $\F$ is the automaton defined by a covering $f:\M_1\arr\M$ of a
topological space by its open subset, then $\M_n$ is the domain of
the $n$th iterate $f^n$ of $f$ and the automaton $\F^n$ is defined
by the partial self-covering $f^n:\M_n\arr\M$.
\end{examp}

\begin{examp}
If the topological automaton is the dual Moore diagram of an
automaton, then the topological automaton $\F^n=(\M, \M_n, f^n,
\iota^n)$ is the dual Moore diagram of the same automaton over the
alphabet $\alb^n$. Analogous statement holds for topological
automata describing wreath recursions on groups.
\end{examp}

\subsection{Three inverse limits of a topological automaton}
Let $\F=(\M, \M_1, f, \iota)$ be a topological automaton.
Iterations of $\F$ produce the following infinite commutative
diagram of topological spaces. We consider $\M_n$ as regular
topological spaces even if $\M$ is an orbispace, i.e., consider
the underlying spaces only.
\[\begin{array}{ccccccc}
\ddots & & \vdots & & \vdots & & \vdots\\
&  & \mapdown{f_4} & & \mapdown{f_3} & & \mapdown{f_2}\\
\ldots & \stackrel{\iota_4}{\arr} & \M_4 &
\stackrel{\iota_3}{\arr}& \M_3 &
\stackrel{\iota_2}{\arr} & \M_2\\
& & \mapdown{f_3} & & \mapdown{f_2} & & \mapdown{f_1}\\
\ldots & \stackrel{\iota_3}{\arr} & \M_3 &
\stackrel{\iota_2}{\arr}&
\M_2 & \stackrel{\iota_1}{\arr} & \M_1\\
& & \mapdown{f_2} & & \mapdown{f_1} &  & \mapdown{f}\\
\ldots & \stackrel{\iota_2}{\arr} & \M_2 &
\stackrel{\iota_1}{\arr}&\M_1&\stackrel{\iota}{\arr}&\M\end{array}\]

Denote by $\lim_f\F$ the inverse limit of the columns of this
diagram. The limit obviously does not depend on the choice of the
column and the maps $\iota_n$ between the columns induce a
continuous map $\iota_\infty:\lim_f\F\arr\lim_f\F$. Similarly,
denote by $\lim_\iota\F$ the inverse limit of the rows. The maps
$f_n$ induce then a continuous map
$f_\infty:\lim_\iota\F\arr\lim_\iota\F$, which is a covering in
the regular case.

We may also consider the inverse limit of the whole diagram, which
we will denote $\lim_{f, \iota}\F$. The ``diagonal'' identical map
between the corners $\M_n$ of the commutative squares
\[\begin{array}{ccc}\M_{n+1} & \stackrel{\iota_n}{\arr} & \M_n\\
\mapdown{f_n} & \swarrow & \mapdown{f_{n-1}} \\
\M_n & \stackrel{\iota_{n-1}}{\arr} & \M_{n-1}\end{array}\]
induces a homeomorphism $\Delta$ of $\lim_{f, \iota}\F$. The
following is straightforward.

\begin{proposition}
The space $\lim_{f, \iota}\F$ is homeomorphic to the inverse limit
of the sequence
\[\ldots\stackrel{f_\infty}{\longleftarrow}\lim_\iota\F
\stackrel{f_\infty}{\longleftarrow}\lim_\iota\F\stackrel{f_\infty}{\longleftarrow}
\lim_\iota\F,\] and to the inverse limit of the sequence
\[\ldots\stackrel{\iota_\infty}{\longleftarrow}\lim_f\F
\stackrel{\iota_\infty}{\longleftarrow}\lim_f\F\stackrel{\iota_\infty}{\longleftarrow}
\lim_f\F.\] The homeomorphism $\Delta$ is induced by the action of
$f_\infty$ on the first inverse limit and the homeomorphism
$\Delta^{-1}$ is induced by the action of $\iota_\infty$ on the
second inverse limit.
\end{proposition}

\begin{defi}
Let $\F=(\M, \M_1, f, \iota)$ be a regular topological automaton.
A \emph{forward $\F$-orbit} is a sequence \(x_1, x_2, \ldots,\) of
points of $\M_1$ such that
\[f(x_n)=\iota(x_{n+1})\]
for all $n=1,2,\ldots$.

A \emph{backward $\F$-orbit} is a sequence \(x_1, x_2, \ldots,\)
of points of $\M_1$ such that
\[f(x_{n+1})=\iota(x_n)\]
for all $n=1,2,\ldots$.

A \emph{bilateral $\F$-orbit} is a sequence $\ldots, x_{-1}, x_0,
x_1, x_2, \ldots$ such that
\[f(x_n)=\iota(x_{n+1})\]
for all $n\in\Z$.
\end{defi}

The choice of the names ``backward'' and ``forward'' is almost
arbitrary. We use the choice given in the definition, since
iteration of partial self-coverings is our main motivation. It is,
however, also natural to use the opposite terminology, like it is
done in~\cite{katsura:generalizing,katsura:examples}, especially
in the setting of automata theory, groupoids or operator algebras.

The spaces of forward, backward and bilateral $\F$-orbits is
endowed with the topology of a subset of the corresponding direct
powers of $\M_1$.

The following description of the inverse limits is a direct
corollary of Proposition~\ref{pr:iteratreg}.

\begin{proposition}
\label{pr:spoforbits} The spaces $\lim_\iota\F$, $\lim_f\F$ and
$\lim_{f, \iota}\F$ are homeomorphic to the spaces of forward,
backward and bilateral $\F$-orbits, respectively. The maps
$f_\infty, \iota_\infty$ and $\Delta$ are induced by the shifts on
the corresponding spaces of orbits.
\end{proposition}

\begin{examp}
If $\M$ is a topological space and $\iota$ is an embedding (i.e.,
if $\F$ is a partial self-covering), then $\lim_\iota\F$ is the
intersection of the domains $\M_n$ of $f^n$.
\end{examp}

\begin{examp}
Let $f\in\C(z)$ be a hyperbolic rational function. Let
$U\subset\CS$ be a closed set such that $U$ does not intersect the
union of the attracting cycles of $f$, $U$ contains the Julia set
of $f$, and $f^{-1}(U)\subseteq U$. Let $\F=(U, f^{-1}(U), f, id)$
be the corresponding topological automaton. Then $\lim_{id}\F$ is
the Julia set of $f$, and $f_\infty$ is the restriction of $f$
onto its Julia set.

The space of backward orbits $\lim_f\F$ for rational functions was
studied (in greater generality)
in~\cite{lyubichminsk,kaimlyubich:laminations}.
\end{examp}

\subsection{Definition of the iterated monodromy group}
\label{ss:defimg} The definition of the iterated monodromy group
of a topological automaton almost coincides with the definition of
the iterated monodromy group of a partial self-covering
(especially in its orbispace version). Here we give a short
overview of the definitions for regular (non-orbispace) case. For
more details and for the definition in the case of orbispace
topological automata, see~\cite{nek:book} (the map $\iota$ is
considered in~\cite{nek:book} to be an embedding of orbispaces,
but this fact is never used).

Let $\F=(\M, \M_1, f, \iota)$ be a topological automaton, and
suppose that $\M$ is path connected and locally path connected.

Choose a basepoint $t\in\M$, and consider the sequence of the
coverings
\[\M\stackrel{f}{\longleftarrow}\M_1\stackrel{f_1}{\longleftarrow}
\M_2\stackrel{f_2}{\longleftarrow}\cdots,\] and denote $f^n=f\circ
f_1\circ\cdots\circ f_{n-1}$, and $f^{-n}=(f^n)^{-1}$. The
fundamental group $\pi_1(\M, t)$ acts on each of the sets
$f^{-n}(t)\subset\M_n$ by the monodromy action: the image of a
point $z\in f^{-n}(t)$ under the action of a loop
$\gamma\in\pi_1(\M, t)$ is the endpoint of the unique lift of
$\gamma$ by $f^n$ that starts at $z$.

The union $T=\bigcup_{n\ge 0}f^{-n}(t)$ is called the
\emph{preimage tree} of the point $t$. We define vertex adjacency
in $T$ in the natural way, so that a point $z\in f^{-n}(t)$ is
connected by an edge to the point $f_{n-1}(z)\in f^{-(n-1)}(t)$.
It is easy to see that the action of the fundamental group on the
sets $f^{-n}(t)$ is an action by automorphisms of the preimage
tree. This action is called the \emph{iterated monodromy action}.

\begin{defi}
The \emph{iterated monodromy group} of a topological automaton
$\F$ is the quotient of the fundamental group of $\M$ by the
kernel of its action on the tree of preimages.
\end{defi}

\subsection{Coding tree}\label{ss:codingtree}
Exactly as in the case of partial self-coverings, the iterated
monodromy group of a topological automaton can be computed using a
natural self-similarity structure on it. We repeat here the
constructions of~\cite[Sections~5.1--2]{nek:book} simplified to
the case of regular topological automata.

For a covering $P:\mathcal{X}_1\arr\mathcal{X}$, a path $\gamma$
in $\mathcal{X}$ and a preimage $z\in\mathcal{X}_1$ of the
beginning of $\gamma$, we denote by $P^{-1}(\gamma)_z$ the lift of
$\gamma$ by $P$ starting at $z$.

Let $\F=(\M, \M_1, f, \iota)$ be a topological automaton with a
path connected and locally path connected base space $\M$. Let
$\alb$ be an alphabet of size equal to the degree of the covering
$f$. Choose a basepoint $t\in\M$ and a bijection
$\Lambda_1:\alb\arr f^{-1}(t)$. Also choose for every $x\in\alb$ a
path $\ell_x$ in $\M$ starting at $t$ and ending in
$\iota(\Lambda_1(x))$.

Define now inductively the points $\Lambda_1(v)\in\M_1$, curves
$\ell^1_v$ in $\M_1$ for $v\in\xs\setminus(\alb^1\cup\alb^0)$, and
$\ell_v$ in $\M$ for $v\in\xs$, by the rules
\begin{equation}
\label{eq:codingtree}
\ell^1_{xvy}=f^{-1}(\ell_{xv})_{\Lambda_1(vy)},\quad
\ell_v=\iota(\ell^1_v),\quad \text{$\Lambda_1(v)$ is the end of
$\ell^1_v$}.
\end{equation} In other words, we lift the curves
of $\M$ by $f$ to curves in $\M_1$ and then push them back into
$\M$ by $\iota$. In this way we get a tree of curves $\ell_v$ in
$\M$ with the root in $t$ and $d$ trees of curves $\ell^1_v$ with
the roots in $f^{-1}(t)$.

The curve $\ell^1_{xv}$ connects $\Lambda_1(v)$ with
$\Lambda_1(xv)$; the curve $\ell_{xv}$ connects the point
$\iota(\Lambda_1(v))$ with the point $\iota(\Lambda_1(xv))$. It
follows from the definition that \[f(\Lambda_1(x_1\ldots
x_k))=\iota(\Lambda_1(x_1\ldots x_{k-1})),\] since the curve
$\ell^1_{x_1\ldots x_k}$ ends in $\Lambda_1(x_1\ldots x_k)$, while
its $f$-image
\[f(\ell^1_{x_1\ldots x_k})=\ell_{x_1\ldots
x_{k-1}}=\iota(\ell^1_{x_1\ldots x_{k-1}})\] ends in
$\iota(\Lambda(x_1\ldots x_{k-1}))$. Consequently, by
Proposition~\ref{pr:iteratreg}, the sequence
\[(\Lambda_1(x_1), \Lambda_1(x_1x_2), \ldots, \Lambda_1(x_1x_2\ldots x_{n-1}),
\Lambda_1(x_1x_2\ldots x_n))\] defines a point of $\M_n$, which we
will denote by $\Lambda(x_1x_2\ldots x_n)$.

\begin{proposition}
The map $\Lambda:\xs\arr T$ is an isomorphism of the tree of words
with the preimage tree.
\end{proposition}

We call the isomorphism $\Lambda$ the \emph{coding} of the
preimage tree, defined by the connecting paths $\ell_x$.

\begin{proof}
A direct corollary of the construction and
Proposition~\ref{pr:iteratreg}.
\end{proof}

\subsection{Computation of the iterated monodromy group}
\begin{theorem}
\label{th:staction} Let $\Lambda:\xs\arr T$ be the coding defined
by a collection of paths $\ell_x$ connecting the basepoint $t$ to
$\iota(\Lambda(x))$. Then the action of $\pi_1(\M, t)$ on $\xs$,
obtained by conjugation of the iterated monodromy action by the
isomorphism $\Lambda$, is defined by the following recurrent rule:
\[\gamma(xv)=y\left(\ell_y^{-1}\iota(f^{-1}(\gamma)_{\Lambda(x)})\ell_x\right)(v),\]
where $y=\gamma(x)$ is the end of $f^{-1}(\gamma)_{\Lambda(x)}$.
\end{theorem}

\begin{remark}
Here and throughout the paper we multiply paths as functions: in a
product $\gamma_1\cdot\gamma_2$ the path $\gamma_2$ is passed
before $\gamma_1$.
\end{remark}

The proof of the above theorem is the same as in the case of
partial self-coverings, see~\cite[Proposition~5.2.2]{nek:book}.

\begin{defi}
The self-similar action of the iterated monodromy group on $\xs$
described in Theorem~\ref{th:staction} is called the
\emph{standard action} (\emph{defined by the bijection
$\Lambda_1:\alb\arr f^{-1}(t)$ and connecting paths $\ell_x$}).
\end{defi}

\begin{proposition}
\label{pr:stactunique} The standard action of the iterated
monodromy group does not depend on the choice of the bijection
$\Lambda_1$ and the connecting paths $\ell_x$, up to equivalence
of self-similar groups. Any self-similar action equivalent to a
standard action is a standard action.
\end{proposition}

In other words, the iterated monodromy group $\img{\F}$ has a
natural well defined self-similarity structure.

\begin{proof}
A change of the bijection is equivalent to post-conjugation of the
wreath recursion by an element of $\symm$. Changing the set of
connecting paths $(\ell_x)_{x\in\alb}$ to a set
$(\ell_x')_{x\in\alb}$ (for a fixed bijection $\Lambda_1$)
corresponds, by Theorem~\ref{th:staction}, to post-conjugating the
wreath recursion by the element
$(\ell_x^{-1}\ell_x')_{x\in\alb}\in(\pi_1(\M, t)^\alb$.
\end{proof}

Standard actions of the iterated monodromy groups can be also
defined using the associated virtual endomorphisms.

Let $(\M, \M_1, f, \iota)$ be, as before, a topological automaton
with a path connected and locally path connected space $\M$. We
will assume now that $\M_1$ is also path connected. Fix some
basepoint $t\in\M$ and a point $t_1\in f^{-1}(t)$. Choose a path
$\ell$ in $\M$ from $t$ to $\iota(t_1)$. Let $G_1$ be the subgroup
of $\pi_1(\M, t)$ of loops $\gamma$ such that the lift
$f^{-1}(\gamma)_{t_1}$ is also a loop. The subgroup $G_1$ is of
index $d$ in $\pi_1(\M, t)$ (and is isomorphic to $\pi_1(\M_1)$).

\begin{defi}
The \emph{virtual endomorphism} of $\pi_1(\M, t)$,
\emph{associated} with the topological automaton is the
homomorphism
\[\phi:G_1\arr\pi_1(\M, t):\gamma\mapsto\ell^{-1}\iota(f^{-1}(\gamma)_{t_1})\ell.\]
\end{defi}
It is easy to check that the associated endomorphism does not
depend, up to conjugacy of virtual endomorphisms, on the choice of
the preimage $t_1$ and of the connecting path $\ell$. Moreover, it
does not depend on the choice of the basepoint $t$, if we identify
the fundamental groups with different basepoints in the standard
way, using connecting paths.

Similarly to the case of partial self-coverings
(see~\cite[Proposition~5.1.2]{nek:book}), one can show that the
standard action of the iterated monodromy group of a topological
automaton is equivalent to the self-similar action defined by the
associated virtual endomorphism.

\subsection{Combinatorial equivalence}
\label{ss:combequiv}

\begin{defi}
We say that two topological automata with path-connected base
(orbi)spaces are \emph{combinatorially equivalent} if their
iterated monodromy groups are equivalent as self-similar groups.
\end{defi}

Here we consider faithful iterated monodromy groups, and not just
the self-similarity on the fundamental group.

\begin{proposition}
\label{pr:combequiv} Let $\F=(\M, \M_1, f, \iota)$ and
$\mathcal{F'}=(\M', \M_1', f', \iota')$ be topological automata
with path-connected base spaces $\M$ and $\M'$. Let
$\phi:\M'\arr\M$ be a continuous map and suppose that the covering
$f':\M_1'\arr\M'$ is the pullback of $f$ by $\phi$. Let then
$\phi_1:\M'_1\arr\M_1$ be the map making the diagram
\[\begin{array}{ccc}\M_1' & \stackrel{\phi_1}{\arr} & \M_1\\
\mapdown{f'} & & \mapdown{f}\\
\M' & \stackrel{\phi}{\arr} & \M\end{array}\] commutative. Suppose
that the diagram
\begin{equation}\label{eq:fundgr}\begin{array}{ccc}\pi_1(\M_1') &
\stackrel{(\phi_1)_*}{\arr} & \pi_1(\M_1)\\
\mapdown{\iota'_*} & & \mapdown{\iota_*}\\
\pi_1(\M') & \stackrel{\phi_*}{\arr} &
\pi_1(\M)\end{array}\end{equation} is also commutative up to an
inner automorphism of $\pi_1(\M)$, and the map $\phi_*$ is an
epimorphism. Then the topological automata $\F$ and $\F'$ are
combinatorially equivalent.

In particular, homotopically equivalent topological automata are
combinatorially equivalent.
\end{proposition}

\begin{proof}
Let $t'\in\M'$ and $t\in\M$ be such that $t=\phi(t')$. Choose a
bijection $\Lambda_1':\alb\arr (f')^{-1}(t')$ and a collection of
connecting paths $\ell_x$, defining the standard action of the
iterated monodromy group of $\F'$. Then $\phi_1\circ\Lambda_1'$
and $\phi(\ell_x)$ is a bijection and a collection of connecting
paths, defining some standard action of $\img{\F}$. Both standard
actions are unique up to an equivalence of self-similar actions.
Let $\psi':\pi_1(\M')\arr\pi_1(\M')^\alb\rtimes\symm$ and
$\psi:\pi_1(\M)\arr\pi_1(\M)^\alb\rtimes\symm$ be the associated
wreath recursions.

It follows then from commutativity of the
diagram~\eqref{eq:fundgr} and Theorem~\ref{th:staction} that if
$\phi_*(g')=g$, then $\psi(g)$ is obtained from $\psi'(g')$ by
applying $\phi_*$ and a fixed inner automorphism of $\pi_1(\M)$ to
every coordinate of $\pi_1(\M')^\alb$. But this implies that the
iterated monodromy groups of $\F'$ and $\F$ are equivalent.
\end{proof}

\subsubsection{Moore diagrams of the standard action}
The process of finding the standard self-similarity on the
iterated monodromy group $\img{\F}$ is naturally interpreted,
using Proposition~\ref{pr:combequiv}, as passing to a topological
automaton that is a Moore diagram of a self-similar group and is
combinatorially equivalent to $\F$.

Let $\F=(\M, \M_1, f, \iota)$ be a topological automaton such that
$\M$ is path connected and locally simply connected. Let
$S=\{\gamma_i\}_{i\in I}$ be a generating set of the fundamental
group $\pi_1(\M, t)$. Let $\Gamma$ be a rose of loops $g_i$ with a
basepoint $t_0$, for $i\in I$ and let $\phi:\Gamma\arr\M$ be a map
such that $\phi(g_i)=\gamma_i$, $\phi(t_0)=t$. Then the map
$\phi:\Gamma\arr\M$ induces a surjective map of the fundamental
groups.

Let the covering $f':\Gamma_1\arr\Gamma$ and the map
$\phi_1:\Gamma_1\arr\M_1$ be obtained by taking pullback of the
covering $f$ by the map $\phi$. If we find a map
$\iota':\Gamma_1\arr\Gamma$ making the diagram
\begin{equation}\label{eq:pi1diagram}
\begin{array}{ccc}\pi_1(\Gamma_1) & \stackrel{{\phi_1}_*}{\arr} & \pi_1(\M_1)\\
\mapdown{\iota'_*} & & \mapdown{\iota_*}\\
\pi_1(\Gamma) & \stackrel{\phi_*}{\arr} &
\pi_1(\M)\end{array}
\end{equation} commutative, then the
one-dimensional topological automaton $(\Gamma, \Gamma_1, f',
\iota')$ will be combinatorially equivalent to $\F$, by
Propositions~\ref{pr:combequiv}.

If $\iota'$ is such that $\iota'((f')^{-1}(t_0))=\{t_0\}$, then
the automaton $(\Gamma, \Gamma_1, f', \iota')$ is the dual Moore
diagram of a wreath recursion that defines a standard action of
$\img{\F}$, by Propositions~\ref{pr:combequiv}
and~\ref{pr:stactunique}. Conversely, any dual Moore diagram
associated with the wreath recursion of a standard action of
$\img{\F}$ can be obtained in this way. If the standard action is
defined by connecting paths $\ell_x$, then for a lift
$h_z\subset\Gamma_1$ of a loop $g_i$ of $\Gamma$ we define
$\iota'(h_z)$ to be a loop $g$ such that
$\phi(g)=\ell_x^{-1}\iota(\phi_1(h_z))\ell_y$, where $y$ is the
beginning and $x$ is the end of $\phi_1(h_z)$. It is easy to check
that so defined map $\iota'$ makes the
diagram~\eqref{eq:pi1diagram} commutative up to an inner
automorphism of $\pi_1(\M)$.

Consequently, Proposition~\ref{pr:combequiv} is a complete
description of combinatorial equivalence. Two automata are
combinatorially equivalent if and only if there exists a third
automaton, combinatorial equivalence of which to the first two can
be established using Proposition~\ref{pr:combequiv}.

\section{Contracting automata}
\label{s:contracting}
\subsection{Self-similar $G$-spaces}\label{ss:ssgspaces}
Let us redefine the notion of a topological automaton in terms of
actions of groups on topological spaces. This approach will help
us to use the techniques of self-similar groups, and will include
orbispace automata into our consideration without heavy use of the
theory of orbispaces.

Let $\bim$ be a covering bimodule over a group $G$ (see
Definition~\ref{def:covbim}) and let $\X$ be a topological space
with a right action of $G$ by homeomorphisms. Then the tensor
product $\X\otimes\bim$ is defined as the quotient of the direct
product $\X\times\bim$ of topological spaces (where $\bim$ is
discrete) by the identifications
\[\xi\cdot g\otimes x=\xi\otimes g\cdot x.\]
The space $\X\times\bim$ is a right $G$-space with respect to the
action
\[(\xi\otimes x)\cdot g=\xi\otimes(x\cdot g).\]

\begin{defi}
A right $G$-space $\X$ is said to be \emph{$\bim$-invariant} (or
\emph{self-similar}) if the right $G$-spaces $\X$ and
$\X\otimes\bim$ are conjugate, i.e., if there exists a
$G$-equivariant homeomorphism $I:\X\otimes\bim\arr\X$. It is
called \emph{$\bim$-semi-invariant} if there exists a
$G$-equivariant continuous map $I:\X\otimes\bim\arr\X$
\end{defi}

An example of a self-similar $G$-space for a contracting group $G$
is the limit $G$-space $\limg$, where the conjugacy $I$ maps
$\xi\otimes x$, for $\xi\in\limg$ and $x\in\alb$, to the point of
$\limg$ represented by $\ldots x_2x_1g(x)\cdot g|_x$, if $\xi$ is
represented by $\ldots x_2x_1\cdot g$
(see~\cite[Section~3.4]{nek:book} and
Subsection~\ref{ss:contractinglimsp} of our paper).

\begin{lemma}
\label{l:propermap} Let $\X_1, \X_2$ be locally compact,
Hausdorff, proper, and co-compact right $G$-spaces. Then every
$G$-equivariant map $\Phi:\X_1\arr\X_2$ is proper, i.e.,
$\Phi^{-1}(C)$ is compact for every compact $C\subset\X_2$.
\end{lemma}

Recall that an action of $G$ on $\X$ is said to be proper if for
every compact subset $C\subset\X$ the set of elements $g\in G$
such that $C\cdot g\cap C\ne\emptyset$ is finite. It is called
co-compact if there exists a compact set $K$ intersecting every
$G$-orbit.

\begin{proof}
Let $K\subset\X_1$ be a compact set such that $\X_1=\bigcup_{g\in
G}K\cdot g$. Let $C\subset\X_2$ be any compact set. The set
$A=\{g\in G\;:\;\Phi(K)\cdot g\cap C\ne\emptyset\}$ is finite by
compactness of $\Phi(K)\cup C$ and properness of the action of $G$
on $\X_2$. Then
\[\Phi^{-1}(C)\subseteq\bigcup_{g\in A}K\cdot g,\]
hence $\Phi^{-1}(C)$ is compact.
\end{proof}

\begin{lemma}\label{l:properproduct}
Let $\X$ be a locally compact Hausdorff right $G$-space and let
$\bim$ be a covering $G$-bimodule.

If the action of $G$ on $\X$ is proper and co-compact, then the
action of $G$ on $\X\otimes\bim$ is also proper and co-compact.
\end{lemma}

\begin{proof}
Let $K\subset\X$ be an open set with compact closure such that
$\X=\bigcup_{g\in G}K\cdot g$. Then every point of $\X\otimes\bim$
can be written in the form $\xi\otimes x$ for $\xi\in K$ and
$x\in\bim$. Let $\alb$ be a basis of $\bim$. Then
$\X\otimes\bim=\bigcup_{g\in G}(K\otimes\alb)\cdot g$. The set
$K\otimes\alb$ is compact, hence the action of $G$ on
$\X\otimes\bim$ is co-compact.

For every $x\in\bim$ the set $K\otimes x\subset\X\otimes\bim$ is
open, since its preimage \[\bigcup_{g\in G}(K\cdot g, g^{-1}\cdot
x)\] in $\X\times\bim$ is open.

Let $C\subset\X\otimes\bim$ be a compact set. There exists a
finite set $M\subset\bim$ such that $C\subset\bigcup_{x\in M}
K\otimes x$. Suppose that $\xi_1, \xi_2\in C$ and $g\in G$ are
such that $\xi_1\cdot g=\xi_2$. Then there exist $\zeta_1,
\zeta_2\in K$ and $x_1, x_2\in M$ such that $\xi_i=\zeta_i\otimes
x_i$ for $i=1,2$. Then $\zeta_1\otimes x_1\cdot g=\zeta_2\otimes
x_2$, which means that there exists $h\in G$ such that
\[\zeta_1=\zeta_2\cdot h, \quad h\cdot x_1=x_2\cdot g^{-1}.\] The first
equality and properness of the action of $G$ on $\X$ implies that
the set of possible $h$ is finite. Then the second equality and
freeness of the right action on $\bim$ implies that the set of
possible values of $g$ is also finite.
\end{proof}

Let $\X$ be a right $G$-space, and let $I:\X\otimes\bim\arr\X$ be
an equivariant continuous map. If the action of $G$ on $\X$ is
proper, then the associated groupoid of germs of the action of $G$
on $\X$ is an atlas of some orbispace $\M$.

Fix a basis $\alb$ of the bimodule $\bim$ (see
Definition~\ref{def:covbim}). Then we have the associated action
of $G$ on $\alb$, hence we get a covering of the orbispace $\M$ by
the orbispace $\M_1$ of the action
\[g(\xi, x)=(\xi\cdot g^{-1}, g(x))\]
of $G$ on $\X\times\alb$. The covering map $p:\M_1\arr\M$ is
induced by the projection map $P:\X\times\alb\arr\X$.

The semi-conjugacy $I:\X\otimes\bim\arr\X$ naturally induces a
functor of the groupoids of germs, hence it defines a morphism
$\iota:\M_1\arr\M$ of the orbispaces. More explicitly, the functor
maps the germ of the action of $g\in G$ at a point $(\xi, x)$ to
the germ of the action of $g|_x$ at $I(\xi\otimes x)$.

\begin{defi}
The constructed automaton $(\M, \M_1, p, \iota)$ is the
\emph{automaton associated with the $G$-space $\X$ and the
semiconjugacy $I:\X\otimes\bim\arr\X$}.
\end{defi}

\subsection{Self-similar $G$-spaces from topological automata}
\label{ss:ssGtopaut}

Suppose that $\F=(\M, \M_1, p, \iota)$ is a topological automaton
such that the space $\M$ is compact, path connected and
semi-locally simply connected (resp.\ developable, if it is an
orbispace). Recall that a topological space $\M$ is semi-locally
simply connected, if for every point $x\in\M$ there exists a
neighborhood $U$ such that every loop in $U$ is homotopic in $\M$
to a point.

The universal covering $\wt{\M}$ of $\M$ is defined as the space
of homotopy classes of paths starting at a fixed basepoint $t$.
The fundamental group $\pi_1(\M, t)$ acts on $\wt{\M}$ in the
usual way: by appending loops to the paths. The action is
co-compact if $\M$ is compact. It is proper by semi-local simple
connectedness and local compactness of $\M$.

The associated bimodule $\bim_{\F}$ over $\pi_1(\M, t)$ is the set
of pairs $(\ell, z)$, where $z\in p^{-1}(t)$ and $\ell$ is a
homotopy class of a path starting in $t$ and ending in $\iota(z)$
(see~\cite[Section~5.1.4]{nek:book}). The fundamental group
$\pi_1(\M, t)$ acts on $\bim_{\F}$ on the right by appending paths
\[(\ell, z)\cdot\gamma=(\ell\gamma, z),\]
and on the left by taking lifts by $p$:
\[\gamma\cdot(\ell, z)=(\iota(p^{-1}(\gamma)_z)\ell, \gamma(z)),\]
where $\gamma(z)$ is the end of $p^{-1}(\gamma)_z$ (i.e., the
image of $z$ under the action of $\gamma$). Recall that in a
product of paths $\ell\gamma$ the path $\gamma$ is passed before
$\ell$.

If $\xi\in\wt{\M}$ is a point represented by a path $\alpha$
starting at $t$, and $(\ell, z)$ is an element of $\bim_{\F}$,
then define $I(\xi\otimes(\ell, z))$ to be the point of $\wt{\M}$
represented by the path $\iota(p^{-1}(\alpha)_z)\ell$.

\begin{proposition}
The map $I(\xi\otimes(\ell, z))=\iota(p^{-1}(\alpha)_z)\ell$ is a
well defined $\pi_1(\M, t)$-equivariant continuous map from
$\wt{\M}\otimes\bim_{\F}$ to $\wt{\M}$.
\end{proposition}

\begin{proof}
Equivariance and the fact that $I$ is well defined follows
directly from the definitions of the actions of $\pi_1(\M, t)$ on
$\wt{\M}$ and $\bim_{\F}$. Continuity follows from continuity of
the map $\iota$ and branches of $p^{-1}$.
\end{proof}

\begin{proposition}
The automaton $\F=(\M, \M_1, p, \iota)$ is isomorphic to the
automaton associated with the $\pi_1(\M, t)$-space $\wt{\M}$ and
the equivariant map $I:\wt{\M}\otimes\bim_{\F}\arr\wt{\M}$.
\end{proposition}

Here two automata $\F=(\M, \M_1, p, \iota)$ and $\F'=(\M', \M_1',
p', \iota')$ are called isomorphic if there exist homeomorphisms
$\psi:\M\arr\M'$ and $\psi_1:\M_1\arr\M_1'$ such that
$p'\circ\psi_1=\psi\circ p$ and
$\iota'\circ\psi_1=\psi\circ\iota$.

\begin{proof}
Fix a basis $\alb=\{x_z=(\ell_z, z)\;:\;z\in p^{-1}(t)\}$ of
$\bim_{\F}$. Let $(\xi, x_z)$ be a point of $\wt{\M}\times\alb$
and suppose that $\xi$ is represented by a path $\alpha$. Define
$\Psi_1(\xi, x_z)\in\M_1$ to be the end of the path
$p^{-1}(\alpha)_z$. For every $\gamma\in\pi_1(\M, t)$ we have
\[\Psi_1(\alpha\cdot\gamma^{-1}, \gamma(z))=\Psi_1(\alpha, z),\] since
\[p^{-1}(\alpha\cdot\gamma^{-1})_{\gamma(z)}p^{-1}(\gamma)_z=
p^{-1}(\alpha\cdot\gamma^{-1}\gamma)_z.\]

In the other direction, suppose that $\Psi_1(\xi_1,
x_{z_1})=\Psi_1(\xi_2, x_{z_2})$ and $\xi_1, \xi_2$ are
represented by paths $\alpha_1$ and $\alpha_2$. Then the endpoints
of the paths $p^{-1}(\alpha_1)_{z_1}$ and $p^{-1}(\alpha_2)_{z_2}$
coincide, hence the path
\[\gamma=p((p^{-1}(\alpha_1)_{z_1})^{-1}p^{-1}(\alpha_2)_{z_2})=\alpha_1^{-1}\alpha_2\]
is an element of $\pi_1(\M, t)$. Then $\xi_2=\xi_1\cdot\gamma$ and
$\gamma^{-1}(z_1)=z_2$, since the path
\[(p^{-1}(\alpha_2)_{z_2})^{-1}p^{-1}(\alpha_1)_{z_1}\] is a lift of
$\gamma^{-1}$ by $p$.

It follows that $\Psi_1$ induces a homeomorphism $\psi_1$ between
the quotient of $\wt{\M}\times\alb$ by the action $(\xi,
x_z)\mapsto (\xi\cdot\gamma, x_{\gamma^{-1}(z)})$ and $\M_1$. It
is checked now directly that this homeomorphism together with the
natural homeomorphism $\psi:\wt{\M}/\pi_1(\M, t)\arr\M$ satisfies
the definition of an isomorphism of automata.
\end{proof}

Consequently, we will not loose any automaton with connected and
semi-locally simply connected base space $\M$, if we pass to
$G$-spaces and equivariant maps. On the other hand, the
$\bim$-semi-invariant $G$-space $\X$ does not have to be
semi-locally simply connected, thus we can use theory of
$\bim$-semi-invariant $G$-spaces for automata with more general
base spaces $\M$.

\subsection{Iteration of automata associated with $G$-spaces}

Let us describe how automata associated with $\bim$-semi-invariant
spaces are iterated.

The proof of the following lemma follows directly from the
definition of tensor products.

\begin{lemma}
\label{lem:Phinequivariant} Suppose that $I:\X\otimes\bim\arr\X$
is a $G$-equivariant continuous map. Then for every $n\ge 1$ the
map $I^{(n)}:\X\otimes\bim^{\otimes n}\arr\X$ given by
\[I^{(n)}(\xi\otimes x_1\otimes x_2\otimes\cdots\otimes x_n)=
I(\ldots I(I(\xi\otimes x_1)\otimes x_2)\ldots\otimes x_n)\] is
also $G$-equivariant.
\end{lemma}

Let $\X$ be a proper co-compact right $G$-space and let
$I:\X\otimes\bim\arr\X$ be a $G$-equivariant map.

Denote by $\overline{\M_n}$ the orbispace $\X\otimes\bim^{\otimes
n}/G$. Denote by $\M_n$ the orbispace of the action of $G$ on
$\X\times\alb^n$ given by
\[g\cdot (\xi, v)=(\xi\cdot g^{-1}, g(v)).\]

\begin{proposition}
\label{pr:crossotimes} The map $(\xi, v)\mapsto \xi\otimes v$
induces a homeomorphism of the underlying space of $\M_n$ with the
underlying space of $\overline{\M_n}$.
\end{proposition}

\begin{proof}
If the points $(\xi_1, v_1)$ and $(\xi_2, v_2)$ belong to one
orbit of the atlas of $\M_n$ then there exists $g\in G$ such that
$(\xi_1, v_1)=(\xi_2\cdot g^{-1}, g(v_2))$. Then
\[\xi_2\otimes v_2=\xi_2\cdot g^{-1}\otimes g\cdot v_2=\xi_2\cdot
g^{-1}\otimes g(v_2)\cdot g|_{v_2}=\xi_1\otimes v_1\cdot
g|_{v_2},\] i.e., $\xi_1\otimes v_1$ and $\xi_2\otimes v_2$ belong
to one $G$-orbit.

On the other hand, if there exists $h\in G$ such that
$\xi_2\otimes v_2=\xi_1\otimes v_1\cdot h$, then there exists $g$
such that $\xi_2=\xi_1\cdot g$ and $g\cdot v_2=v_1\cdot h$, by the
definition of a tensor product. Then $v_1=g(v_2)$ and
\[(\xi_1, v_1)=(\xi_2\cdot g^{-1}, g(v_2)),\]
i.e., $(\xi_1, v_1)$ and $(\xi_2, v_2)$ belong to one orbit of the
atlas of $\M_n$.
\end{proof}

Even though the underlying spaces of $\M_n$ and $\overline{\M_n}$
are homeomorphic, the orbispaces might be different. The isotropy
groups of $\overline{\M_n}$ are quotients of the corresponding
isotropy groups of $\M_n$.

The following proposition follow directly from the definitions,
see also Proposition~\ref{pr:iteratreg}.

\begin{proposition}
\label{pr:iteratebim} Let $\F=(\M, \M_1, p, \iota)$ be the
automaton associated with a proper right $G$-space $\X$ and a
$G$-equivariant map $I:\X\otimes\bim\arr\X$.

Then the $n$th iteration of $\F$ is the automaton \[\F^{\circ
n}=(\M, \M_n, p_0\circ\cdots\circ p_n,
\iota_0\circ\cdots\circ\iota_n)\] associated with the space $\X$
and the semiconjugacy $I^{(n)}:\X\otimes\bim^{\otimes n}\arr\X$.

The covering $p_n:\M_{n+1}\arr\M_n$ is induced by the
correspondence
\[\xi\otimes v\otimes x\mapsto\xi\otimes v,\]
for $v\in\bim^{\otimes n}$ and $x\in\bim$.

The map $\iota_n:\M_{n+1}\arr\M_n$ is induced by the map
$I_n:\X\otimes\bim^{\otimes (n+1)}\arr\X\otimes\bim^{\otimes n}$
given by
\[I_n(\xi\otimes x\otimes v)=I(\xi\otimes x)\otimes
v\] for $x\in\bim$ and $v\in\bim^{\otimes n}$.
\end{proposition}

The proof of the following proposition is the same as the proof of
Theorem~5.3.1 in~\cite{nek:book}.

\begin{proposition}
Suppose that $\X$ is a path-connected proper right $G$-space and
suppose that $I:\X\otimes\bim\arr\X$ is a $G$-equivariant
continuous map. Then the iterated monodromy group of the
associated topological automaton coincides with the faithful
quotient of the self-similar group $G$.
\end{proposition}

\subsection{Contracting self-similarities}
\begin{defi}
Suppose that $\X$ is a metric space and $G$ acts on it by
isometries, so that the action is proper and co-compact. We say
that an equivariant map $I:\X\otimes\bim\arr\X$ is
\emph{contracting} if there exist $n$ and $0<\lambda<1$ such that
\[d(I^{(n)}(\xi_1\otimes v), I^{(n)}(\xi_2\otimes v))\le\lambda d(\xi_1,
\xi_2)\] for all $\xi_1, \xi_2\in\X$ and $v\in\bim^{\otimes n}$.
\end{defi}

The maps $I^{(n)}$ are defined in Lemma~\ref{lem:Phinequivariant}.

\begin{theorem}
\label{th:approximationlimg} Let $(G, \alb)$ be a contracting
group and let $\bim$ be the associated permutational $G$-bimodule.
Suppose that $\X$ is a locally compact metric space with a
co-compact proper right $G$-action by isometries and let
$I:\X\otimes\bim\arr\X$ be a contracting equivariant map. Then the
projective limit of the $G$-spaces and the $G$-equivariant maps
\[\X\xleftarrow{I_0}\X\otimes\bim\xleftarrow{I_1}\X\otimes\bim^{\otimes 2}
\xleftarrow{I_2}\X\otimes\bim^{\otimes 3}\xleftarrow{I_3}\cdots,\]
is homeomorphic as a $G$-space to the limit $G$-space $\limg$.

For every $x\in\bim$ the maps \[(\xi\otimes v)\mapsto (\xi\otimes
v\otimes x):\X\otimes\bim^{\otimes n}\arr\X\otimes\bim^{\otimes
(n+1)}\] agree with the maps $I_n$ and their limit is the map
$\xi\mapsto\xi\otimes x$ on $\limg$.
\end{theorem}

The maps $I_n:\X\otimes\bim^{\otimes
(n+1)}\arr\X\otimes\bim^{\otimes n}$ were defined in
Proposition~\ref{pr:iteratebim}.

\begin{proof}
Let $K_0\subset\X$ be a compact set such that $\bigcup_{g\in
G}K_0\cdot g=\X$. Choose a basis $\alb$ of $\bim$. There exists a
compact set $K\supseteq K_0$ such that for every $x\in\alb$ and
$\xi\in K$ we have $I(\xi\otimes x)\in K$. One can take, for
instance, the closure of the set of points of the form
$I^{(n)}(\xi\otimes v)$ for $\xi\in K_0$ and $v\in\alb^n$, which
has finite diameter, by contraction of $I$.

Every point of $\X\otimes\bim^n$ can be written in the form
$\xi\otimes v\cdot g$ for $\xi\in K$, $v\in\alb^n$ and $g\in G$.

Hence, every point $\zeta$ of the inverse limit is represented by
a sequence
\[\xi_0\cdot g,\quad\xi_1\otimes x_1\cdot g,
\quad\xi_2\otimes x_2x_1\cdot g,\quad\ldots,\quad \xi_n\otimes
x_n\dots x_2x_1\cdot g,\] for some $g\in G$, $x_i\in\alb$ and
$\xi_n\in K$ such that $I(\xi_n\otimes x_n)=\xi_{n-1}$ for all
$n\ge 1$.

Let us put into correspondence to $\zeta$ the point
$L(\zeta)\in\limg$ represented by the sequence $\ldots x_2x_1\cdot
g$.

Let us show that the map $L$ is well defined. Suppose that we have
the same point $\zeta$ of the inverse limit is represented in two
different ways:
\[\xi_0\cdot g=\eta_0\cdot h,\quad\xi_1\otimes x_1\cdot g=\eta_1\otimes y_1\cdot h,\quad
\xi_2\otimes x_2x_1\cdot g=\eta_2\otimes y_2y_1\cdot
h,\quad\ldots\] for $\xi_n, \eta_n\in K$, $x_n, y_n\in\alb$ and
$g, h\in G$. Then there exists a sequence $g_n\in G$ such that
\[\xi_n=\eta_n\cdot g_n,\qquad g_n\cdot x_n\ldots x_2x_1\cdot
g=y_n\ldots y_2y_1\cdot h\] for all $n\ge 0$. (Notation is
explained after Proposition~\ref{pr:tilesintersection}.) By
compactness of $K$ and properness of the action of $G$ on $\X$,
the set of possible values of the sequence $g_n$ is finite, hence
the sequences $\ldots x_2x_1\cdot g$ and $\ldots y_2y_1\cdot h$
are asymptotically equivalent and represent the same point of
$\limg$.

Let us show that every point of $\limg$ is equal to $L(\zeta)$ for
some point $\zeta$ of the inverse limit of the spaces
$\X\otimes\bim^{\otimes n}$ (i.e., that the map $L$ is onto). Let
$\ldots x_2x_1\cdot g$ be an arbitrary point of the limit
$G$-space $\limg$. For every $n$ chose an arbitrary point $\xi_{n,
n}\in K$ and consider for $k=1, \ldots, n-1$ the points defined
inductively as $\xi_{n, k-1}=I(\xi_{n, k}\otimes x_k)$. Then the
sequence
\[\xi_{n, 0}\cdot g,\quad\xi_{n, 1}\otimes x_1\cdot g,\quad\ldots\quad\xi_{n,
n}\otimes x_n\ldots x_2x_1\cdot g\] agrees with the maps in the
inverse sequence of the spaces $\X\otimes\bim^{\otimes n}$.  We
can choose, by compactness of $K$, an increasing sequence $n_{k,
1}$ such that the sequence $\xi_{n_{k, 1}, 0}$ converges to a
point $\xi_0$. Then we can choose a subsequence $n_{k, 2}$ of
$n_{k, 1}$ such that $\xi_{n_{k, 2}, 1}$ converges to a point
$\xi_1$, etc. In the limit, by continuity of $I$, we get a
sequence
\[\xi_0\cdot g,\quad\xi_1\otimes x_1\cdot g,\quad\ldots,\quad\xi_n\otimes x_n\ldots
x_2x_1\cdot g,\quad\ldots\] representing a point $\zeta$ of the
inverse limit such that $L(\zeta)=\ldots x_2x_1\cdot g$.

Let us show that $L$ is a one-to-one map. Suppose that we have two
sequences
\[\xi_0\cdot g,\quad\xi_1\otimes x_1\cdot g,\quad\xi_2\otimes x_2x_1\cdot
g,\quad\ldots\] and
\[\eta_0\cdot h,\quad\eta_1\otimes y_1\cdot h,\quad\eta_2\otimes y_2y_1\cdot
h,\quad\ldots\] for $\xi_n, \eta_n\in K$, $x_n, y_n\in\alb$ and
$g, h\in G$ such that the sequences $\ldots x_2x_1\cdot g$ and
$\ldots y_2y_1\cdot h$ are asymptotically equivalent. Let $g_n\in
G$ be the sequence implementing the asymptotic equivalence, i.e.,
a sequence with a finite set $A$ of values such that $g_n\cdot
x_n\ldots x_2x_1\cdot g=y_n\ldots y_2y_1\cdot h$. Then the second
point of the inverse limit is written as
\[\eta_0\cdot g_0g,\quad\eta_1\cdot g_1\otimes x_1\cdot g,\quad\eta_2\cdot
g_2\otimes x_2x_1\cdot g,\quad\ldots.\] We have
$I^{(n)}(\xi_n\otimes x_n\ldots x_2x_1\cdot g)=\xi_0\cdot g$ and
$I^{(n)}(\eta_n\cdot g_n\otimes x_n\ldots x_2x_1\cdot
g)=\eta_0\cdot g_0g$. The points $\xi_n$ and $\eta_n\cdot g_n$
belong to the set $K\cup K\cdot A$ of finite diameter. This
implies, by contraction of $I$ that $\xi_0=\eta_0\cdot g_0$. One
proves in the same way that $\xi_n=\eta_n\cdot g_n$ for all $n$,
i.e., that the two points of the inverse limit are the same.

Continuity of the map $L^{-1}$ (i.e., that sequences with long
common beginnings correspond to close points of the inverse limit)
follows directly from the contraction property for $I$.

The maps $I^{(n)}=I_1\circ\cdots\circ I_n:\X\otimes\bim^{\otimes
n}\arr\X$ are proper by Lemmata~\ref{l:properproduct}
and~\ref{l:propermap}, which implies that the inverse limit is
locally compact.

We have constructed an equivariant continuous bijection between
the inverse limit of the spaces $\X\otimes\bim^{\otimes n}$ and
$\limg$. Since both spaces are locally compact and Hausdorff, this
map is a homeomorphism.

The statement about the map $\xi\mapsto\xi\otimes x$ follows
directly from the construction of the homeomorphism $L$.
\end{proof}

We see that if the equivariant map $I:\X\otimes\bim\arr\X$ is
contracting, then the spaces $\X\otimes\bim^{\otimes n}$ are
approximations of the limit $G$-space $\limg$, hence the spaces
$\M_n$ are approximations of the limit space $\lims$ of the group
$G$.

\subsection{Contracting topological automata}

\begin{defi}
\label{def:contractingF} Let $\F=(\M, \M_1, f, \iota)$ be a
topological automaton such that $\M$ is a compact, path connected,
and semi-locally path connected (orbi)space. We say that the
topological automaton $\F$ is \emph{contracting} if there exists a
length structure on $\M$ and $\lambda<1$ such that for every
rectifiable path $\gamma$ in $\M_1$ the path $\iota(\gamma)$ has
length less than $\lambda$ times the length of $\gamma$ (with
respect to the length structure on $\M_1$ equal to the pull back
by $f$ of the length structure on $\M$)
\end{defi}

Let $\F=(\M, \M_1, p, \iota)$ be a contracting automaton. Let
$\wt{\M}$ be the universal covering of $\M$, let $\bim_{\F}$ be
the associated $\pi_1(\M)$-bimodule, and let
$I:\wt{\M}\otimes\bim_{\F}\arr\wt{\M}$ be the self-similarity
defined in Subsection~\ref{ss:ssGtopaut}. Then $\wt{\M}$ has a
natural length structure, which is the lift of the length
structure on $\M$ (length of a curve in $\wt{\M}$ is equal to the
length of its image in $\M$). It follows then from the definition
of the map $I$ and Definition~\ref{def:contractingF} that $I$ is
contracting. This implies that the iterated monodromy group of
$\F$ is a contracting self-similar group. We get now the following
corollary of Theorem~\ref{th:approximationlimg} and
Proposition~\ref{pr:iteratebim}.

\begin{theorem}
\label{th:projlim} Let $\F=(\M, \M_1, f, \iota)$ be a contracting
topological automaton with semi-locally simply connected orbispace
$\M$. Then the iterated monodromy group $\img{\F}$ is contracting
and the system $(\lim_\iota\F, f_\infty)$ is topologically
conjugate to the limit dynamical system $(\lims[\img{\F}], \si)$
of the iterated monodromy group.
\end{theorem}

This means that if a topological automaton $\F=(\M_1, \M, f,
\iota)$ is contracting, then the spaces $\M_n$ can be used as
approximations of the limit space $\lims[\img{\F}]$. The natural
maps $\pi_n:\lim_\iota\F\arr\M_n$ will become more and more
``precise'' in the sense that the difference
$\M_n\setminus\pi_n(\lim_\iota\F)$ and the fibers $\pi_n^{-1}(x)$
become ``smaller''.

\begin{corollary}
\label{cor:imgequiv} If $\F_1$ and $\F_2$ are combinatorially
equivalent contracting automata, then the dynamical systems
$(\lim_\iota\F_1, f_\infty)$ and $(\lim_\iota\F_2, f_\infty)$ are
topologically conjugate.
\end{corollary}

\subsection{Homotopy equivalence and contracting automata}

Contracting topological automata can be simplified using the
following general procedure.

\begin{proposition}
\label{pr:simplification} Let $(G, \alb)$ be a contracting group
and let $\X_1$ and $\X_2$ be metric spaces with proper co-compact
right actions of $G$ by isometries. Suppose that there exists a
contracting $G$-equivariant map $I:\X_1\otimes\bim\arr\X_1$ and
Lipschitz $G$-equivariant maps $F_1:\X_1\arr\X_2$ and
$F_2:\X_2\arr\X_1$. Then there exists $m\ge 1$ and a contracting
$G$-equivariant map $\Psi:\X_2\otimes\bim^{\otimes m}\arr\X_2$.
\end{proposition}

\begin{proof}
For any given $m$ consider the map $\Psi:\X_2\otimes\bim^{\otimes
m}\arr\X_2$ defined by the equality
\[\Psi(\xi\otimes v)=F_1\left(I_1^{(m)}(F_2(\xi)\otimes
  v)\right).\]
It is easy to see that it is $G$-equivariant and well defined (the
latter means that $F_2(\xi)\otimes v$ depends only on $\xi\otimes
v$). If $F_1$ and $F_2$ are Lipschitz with coefficient $L$ and
$I^{(m)}$ is contracting with coefficient $\lambda<L^{-2}$, then
$\Psi$ is contracting.
\end{proof}

\begin{corollary}
\label{cor:hequiv} Let $\F=(\M, \M_1, f, \iota)$ be a contracting
topological automaton such that $\M$ is a finite simplicial
complex (or complex of groups) with a piecewise Riemannian length
structure. Then for any finite simplicial complex (resp.\ complex
of groups) $\M'$ homotopically equivalent to $\M$ there exists $n$
and a contracting automaton $\F'=(\M', \M_1', f', \iota')$
combinatorially equivalent to $\F^{\circ n}$.
\end{corollary}

\begin{proof}
Homotopy equivalence will lift to a pair of
$\pi_1(\M)\cong\pi_1(\M')$-equivariant maps between the universal
coverings of $\M$. By Simplicial Approximation Theorem, we may
assume that these maps are simplicial, hence Lipschitz for some
length structure on $\M'$.
\end{proof}

\section{Simplicial approximations of the limit spaces}
\label{s:simplicialapprox}

\subsection{Topological nucleus}

Let $(G, \alb)$ be a contracting self-similar group with nucleus
$\nuke$. The aim of this section is to find a simple construction
of a proper co-compact $G$-space $\X$ and a contracting
equivariant map $I:\X\otimes\bim\arr\X$.

We assume for simplicity that the group $G$ is finitely generated,
and the action $(G, \alb)$ is self-replicating, i.e., that the
left action of $G$ on the self-similarity bimodule $\bim=\alb\cdot
G$ is transitive.

Defining an equivariant continuous map $I:\X\otimes\bim\arr\X$ is
equivalent to defining a family of continuous maps
$I_x:\X\arr\X:\xi\mapsto I(\xi\otimes x)$ for all $x\in\bim$
satisfying the conditions
\begin{equation}
\label{eq:Phix} I_x(\xi\cdot g)=I_{g\cdot x}(\xi),\qquad I_{x\cdot
  g}(\xi)=I_x(\xi)\cdot g,
\end{equation}
for all $x\in\bim$, $g\in G$ and $\xi\in\X$. The first condition
is equivalent to the condition for the map $I(\xi\otimes
x)=I_x(\xi)$ to be well defined. The second condition is
equivalent to equivariance of $I$.

The iteration $I^{(n)}:\X\otimes\bim^{\otimes n}\arr\X$ is defined
then by the compositions $I_{x_1\otimes\cdots\otimes
  x_n}(\xi)=I_{x_1}\circ\cdots\circ I_{x_n}(\xi)$.

Note that it is enough to define the maps $I_x$ for $x\in\alb$,
since every element of $\bim$ can be written as $x\cdot g$ for
$x\in\alb$ and $g\in G$. Then the only condition to check is
\begin{equation}
\label{eq:phix} I_x(\xi\cdot g)=I_{g(x)}(\xi)\cdot g|_x
\end{equation}
for all $x\in\alb$ and $g\in G$.

Moreover, since we assume that the action is self-replicating,
every element of $\bim$ can be written as $g\cdot x$ for a fixed
$x\in\bim$. Then it is enough to define one map $I_x$ satisfying
condition~\eqref{eq:phix} for all $g$ in the stabilizer of $x$.
This coincides with condition~\eqref{eq:Fmodel} from Introduction.

The simplest example of a $\bim$-semi-invariant $G$-space is the
group $G$ itself with respect to the action by right translations
and the maps
\begin{equation}\label{eq:Phixrestr}I_x(g)=g|_x.\end{equation} Then
condition~\eqref{eq:phix} follows directly from the basic
properties~\eqref{eq:sections} of sections.

Next natural construction will be to choose a finite generating
set $S=S^{-1}$ and consider the corresponding Rips complex, i.e.,
the simplicial complex $\Gamma(G, S)$ with the set of vertices $G$
and the set of simplices equal to the set of subsets $A\subset G$
such that $gh^{-1}\in S$ for all $g, h\in A$. If $S$ is
\emph{self-similar}, i.e., if $S|_x\subset S$ for all $x\in\alb$,
then the maps $I_x$, defined by~\eqref{eq:Phixrestr}, are
simplicial and define a $G$-equivariant map $I:\Gamma(G,
S)\otimes\bim\mapsto\Gamma(G, S)$.

For every finite self-similar set $S$ there exists $n$ such that
$S_n=\bigcup_{v\in\alb^n}S|_v$ is a subset of the nucleus $\nuke$.
Then $I^{(n)}(\Gamma(G, S_n)\otimes\bim^{\otimes n})$ is a
subcomplex of $\Gamma(G, \nuke)$. Consequently, it is sufficient
to consider just the case $S=\nuke$.

Moreover, it may happen that $I(\Gamma(G, \nuke)\otimes\bim)$ is a
proper sub-complex of $\Gamma(G, \nuke)$, and we can then pass to
a smaller complex. Namely, we get a decreasing sequence of
simplicial complexes
\[\Gamma(G, \nuke)\supseteq I(\Gamma(G,
\nuke)\otimes\bim)\supseteq I^{(2)}(\Gamma(G,
\nuke)\otimes\bim^{\otimes 2})\supseteq\cdots,\] which has to
stabilize, since all these complexes are $G$-invariant, and
$\Gamma(G, \nuke)$ is locally finite.

Let us describe the complex $\Gamma=\bigcap_{n\ge
0}I^{(n)}(\Gamma(G, \nuke)\otimes\bim^{\otimes n})$. Since it is
$G$-invariant, it is sufficient to describe the set of simplices
containing the identity.

\begin{proposition}
\label{pr:adjacencyinf} A subset $A\subset\nuke$ containing the
identity element is a simplex of $\Gamma$ if and only if there
exists a sequence $\ldots x_2x_1\in\xmo$ and a sequence
$A_n\subset\nuke$ such that $A_0=A$ and $A_n|_{x_n}=A_{n-1}$ for
all $n\ge 1$.
\end{proposition}

\begin{proof}
A subset $A\subset G$ is a simplex of $\Gamma$ if and only if
there exists a sequence $B_n$ of simplices of $\Gamma(G, \nuke)$
and a sequence of words $v_n\in\alb^n$ such that $B_n|_{v_n}=A$.

If $A_n$ and $\ldots x_2x_1\in\xmo$ satisfy the conditions of the
proposition, then we can take $B_n=A_n$ and $v_n=x_n\ldots
x_2x_1$, and conclude that $A$ is a simplex of $\Gamma$.

Let us prove the other direction of the proposition. Let
$\mathcal{A}_n$ be the set of simplices $A_n$ of $\Gamma(G,
\nuke)$ containing the identity for which there exists $x_n\ldots
x_2x_1\in\alb^n$ such that $A_n|_{x_n\ldots x_2x_1}=A$. Note that
in this case we have $A_n|_{x_n}\in\mathcal{A}_n$, so that we get
a sequence of maps $A_n\mapsto A_n|_{x_n}$ from $\mathcal{A}_n$ to
$\mathcal{A}_{n-1}$. The sets $\mathcal{A}_n$ are finite, and
every element of the inverse limit of the sets $\mathcal{A}_n$
with respect to the described maps gives us sequence $(A_n)_{n\ge
1}$ and $\ldots x_2x_1\in\xmo$ satisfying the conditions of the
proposition. It remains hence to prove that the sets
$\mathcal{A}_n$ are not empty.

If a simplex $B_n$ of $\Gamma(G, \nuke)$ and a word $v\in\alb^n$
are such that $B_n|_v=A$, then there exists $h\in B_n$ such that
$h|_v=1$. Then $B_n\cdot h^{-1}$ is a simplex of $\Gamma(G,
\nuke)$ containing the identity and such that $(B_n\cdot
h^{-1})|_{h(v)}=B_n|_v\cdot (h|_v)^{-1}=A$, i.e., $B_n\cdot
h^{-1}$ is an element of $\mathcal{A}_n$.
\end{proof}

As a direct corollary of Proposition~\ref{pr:adjacencyinf} we get
a more explicit description of the simplices of $\Gamma$.

\begin{corollary}
\label{cor:adjacencygraph} Denote by $\mathcal{B}$ the set of
subsets of $\nuke$ containing the identity. Construct a directed
graph with the set of vertices $\mathcal{B}$ in which a there is
an arrow starting at a vertex $A_1$ and ending in a vertex $A_2$
if there exists $x\in\alb$ such that $A_2=A_1|_x$.

A set $A\in\mathcal{B}$ is a simplex of $\Gamma$ if and only if it
is an end of a directed path starting in a directed cycle in the
constructed graph.
\end{corollary}

\begin{proposition}
\label{pr:Multinucleus} For every finite subset $A\subset G$ there
exists $n$ such that for all words $v\in\xs$ of length at least
$n$ the set $A|_v$ is a simplex of $\Gamma$.
\end{proposition}

\begin{proof}
Replacing $A$ by $A\cdot g^{-1}$ for some $g\in A$, we may assume
that $A$ contains the identity. There exists $n_1$ such that
$A|_v\subset\nuke$, by definition of the nucleus. Then, by
Corollary~\ref{cor:adjacencygraph}, there exists $n_2$ such that
$A|_v|_u$ is a simplex of $\Gamma$ for all words $u$ of length at
least $n_2$. The number $n=n_1+n_2$ satisfies the conditions of
the proposition.
\end{proof}

The simplices of $\Gamma$ have the following geometric
description.

\begin{proposition}
\label{pr:adjacencysets} A subset $A\subset G$ is a simplex of
$\Gamma$ if and only if there exists a point $\xi$ of $\limg$ such
that for every $g\in A$ the point $\xi$ can be represented by
$\ldots x_2x_1\cdot g\in\xmo\times G$.
\end{proposition}

In other words, $A$ is a simplex of $\Gamma$ if and only if the
intersection
\[T_A=\bigcap_{g\in A}\til\cdot g\]
is non-empty.

\begin{proof}
Suppose that for every $g\in A$ there exists a sequence $\ldots
x_2x_1\cdot h\in\xmo\times G$ representing a point of
$\bigcap_{g\in A}\til\cdot g$, i.e., equivalent to some sequences
$w_g=\ldots y_2y_1\cdot g$.

For every $w_g$ there exists then a sequence $h_{n, g}\in\nuke$
such that $h_{n, g}\cdot x_n=y_n\cdot h_{n-1, g}$ and $h_0h=g$.
Denote $A_n=\{h_{n, g}\}_{g\in A}$. Then $A_n|_{x_n}=A_{n-1}$ and
$A_0=A\cdot h^{-1}$, which implies that $A$ is a simplex of
$\Gamma$.

In the other direction, if $A_n$ and $\ldots x_2x_1$ are such that
$A_n|_{x_n}=A_{n-1}$ and $A_0=A\cdot h^{-1}$ for some $h\in G$,
then for every $g\in A$ there exists a sequence $h_{n, g}\in A_n$
such that $h_{n, g}|_{x_n}=h_{n-1, g}$ and $h_{0, g}h=g$. In this
case the sequence $\ldots x_2x_1\cdot h$ is equivalent to the
sequences $\ldots h_{2, g}(x_2)h_{1, g}(x_1)\cdot g$, i.e., the
corresponding point of $\limg$ belongs to every tile $\til\cdot
g$, $g\in A$.
\end{proof}

In view of Proposition~\ref{pr:adjacencysets} we introduce the
following terminology.

\begin{defi}
The complex $\Gamma$ is called the \emph{tiling nerve} of $(G,
\alb)$. The simplices of $\Gamma$ are called \emph{adjacency
subsets} of $G$.
\end{defi}

Denote by $\Gamma_n$ the $G$-space $\Gamma\otimes\bim^{\otimes
n}$. Recall that the maps $I_v:\Gamma\arr\Gamma$, for
$v\in\alb^n$, defining the equivariant maps
$I^{(n)}:\Gamma_n\arr\Gamma$ are simplicial maps given by
\[I_v(g)=g|_v.\]

\begin{defi}
Denote by $J_n(G)$ the quotient $\Gamma_n/G$ (in particular,
$J_0(G)$ is $\Gamma/G$). The topological automaton $(J_0(G),
J_1(G), p, \iota)$ associated with the map
$I:\Gamma\otimes\bim\arr\Gamma$ is called the \emph{topological
nucleus} of the group $(G, \alb)$.
\end{defi}

Recall that the covering $p:J_1(G)\arr J_0(G)$ is induced by the
correspondence $\xi\otimes x\mapsto\xi$ from $\Gamma\otimes\bim$
to $\Gamma$, the map $\iota:J_1(G)\arr J_0(G)$ is induced by the
map $I:\Gamma\otimes\bim\arr\Gamma$.

Note that the restriction of the topological nucleus onto its
one-skeleton coincides with the dual Moore diagram of the nucleus
of $G$. In particular, the one-skeleton of the complex $J_n(G)$ is
the Schreier graph of the action of $G$ on the $n$th level of the
tree $\xs$.

\subsection{Recurrent description of $J_n(G)$}

The spaces $J_n(G)$ can be constructed by the following simple
cut-and-paste procedure.

The barycentric subdivision $\Gamma'$ of $\Gamma$ is isomorphic,
as a simplicial complex, to the realization of the poset (with
respect to inclusion) of the adjacency subsets of $G$.

Let us take as a fundamental domain of the $G$-action on $\Gamma$
the union $T_0$ of the simplices of the barycentric subdivision
$\Gamma'$ containing $1\in G$.

The set of vertices of $T_0$ is the set $\mathcal{A}$ of adjacency
subsets of $G$ containing the identity. The complex $T_0$ is
isomorphic to the geometric realization of the poset $\mathcal{A}$
with respect to inclusion.

For every $g\in\nuke\setminus\{1\}$ consider the subset
$\mathcal{A}_g$ of the poset $\mathcal{A}$ consisting of the
adjacency sets containing $g$. It is a sub-poset of $\mathcal{A}$
and it is equal to the intersection of $\mathcal{A}$ with
$\mathcal{A}\cdot g=\{A\cdot g\;:\;A\in\mathcal{A}\}$. Let $K_{g,
0}$ be the corresponding sub-complex of $T_0$. For every
$A\in\mathcal{A}_g$ the set $A\cdot g^{-1}$ belongs to
$\mathcal{A}_{g^{-1}}$, since $A\cdot g^{-1}\supset\{1, g\}\cdot
g^{-1}=\{g^{-1}, 1\}$.

The map $A\mapsto A\cdot g^{-1}$ is an order-preserving bijection
from $\mathcal{A}_g$ to $\mathcal{A}_{g^{-1}}$. Denote by
$\kappa_{g, 0}:K_{g, 0}\arr K_{g^{-1}, 0}$ the corresponding
isomorphism of the sub-complexes. It coincides with the
restriction of the map $\xi\mapsto\xi\cdot g^{-1}$ onto $T_0\cap
T_0\cdot g=K_{g, 0}$.

It follows then directly from the definitions that the complex
$J_0(G)$ is the quotient of $T_0$ by the identifications
$\kappa_g$ for all $g\in\nuke$.

\begin{proposition}
\label{pr:cutandpaste} Define the complex $T_n$ inductively as the
quotient of $T_{n-1}\times\alb$ by the identifications
\[(\xi, x)\sim (\kappa_{g, n-1}(\xi), g(x))\]
for all $\xi\in K_{g, n-1}$, $x\in\alb$ and
$g\in\nuke\setminus\{1\}$ such that $g|_x=1$.

Define the identification $\kappa_{g, n}$ on the image $K_{g, n}$
of the set \[\bigcup_{x\in\alb, h\in\nuke, h|_x=g}\left(K_{h,
n-1}, x\right)\subset T_{n-1}\times\alb\] in $T_n$ and acting by
the rule
\[\kappa_{g, n}:(\xi, x)\mapsto (\kappa_{h, n-1}(\xi),
h(x)),\] where $h\in\nuke$ is such that $h|_x=g$.

Then the complex $J_n(G)$ is isomorphic to the quotient of $T_n$
by the identifications $\kappa_{g, n}$.

The covering $p_n:J_{n+1}(G)\arr J_n(G)$ is induced by the map
$(\xi, x)\mapsto \xi$ for $\xi\in J_n(G)$ and $x\in\alb$. The map
$\iota_n:J_{n+1}(G)\arr J_n(G)$ is induced by the map $(\xi,
x)\mapsto (\iota_{n-1}(\xi), x)$ for $\xi\in J_n(G)$ and
$x\in\alb$.
\end{proposition}

\begin{proof}
Direct corollary of the definition of the tensor product
$\Gamma\otimes\bim^{\otimes n}$ and
Proposition~\ref{pr:iteratebim}. Here $T_n$ is the fundamental
domain $\bigcup_{v\in\alb^n}T_0\otimes v$ of the action of $G$ on
$\Gamma\otimes\bim^{\otimes n}$.
\end{proof}

The recursive rule described in Proposition~\ref{pr:cutandpaste}
is conveniently encoded by the dual Moore diagram of the nucleus.
Recall that it this diagram the vertices are the letters of the
alphabet $\alb$, and for every $g\in\nuke\setminus\{1\}$ there is
an arrow labeled by $(g, g|_x)$ starting at $x$ and ending in
$g(x)$.

We can interpret now the arrows of the dual Moore diagram of the
nucleus as instructions how to paste together the copies
$(T_{n-1}, x)$ of $T_{n-1}$ into the complex $T_n$, and the
identifications $\kappa_{g, n-1}$ into the identifications
$\kappa_{g, n}$.

Namely, every vertex $x$ corresponds to the piece $(T_{n-1}, x)$.
The arrows labeled by $(g, 1)$ describe
how to paste together the complexes $(T_{n-1}, x)$ into $T_n$: one has to take
the piece corresponding to the beginning of the arrow and attach it by
the map $\kappa_{g, n-1}$ to the piece corresponding to the end of the
arrow.

Every arrow labeled by $(h, g)$ will describe the part of the
identification rule $\kappa_{g, n}$ that maps $(\xi, x)$ to
$(\kappa_{h, n-1}(\xi), y)$, where $x$ is the beginning and $y$ is
the end of the arrow.

\subsection{Topological nucleus as a contracting automaton}
\label{ss:contractingautomaton}

The maps $I^{(n)}:\Gamma_n\arr\Gamma$ are not contracting, since
there always exist simplices $S$ of $\Gamma_n$ mapped
isometrically by $I_v$ for some $v\in\alb^n$. Nevertheless, we can
transform $I^{(n)}$ into a contracting map.

\begin{theorem}
\label{th:contractingautomaton} There exists $n\ge 1$ such that
the map $I^{(n)}:\Gamma\otimes\bim^{\otimes n}\arr\Gamma$ is
$G$-equivariantly homotopic (i.e., is homotopic through
equivariant maps) to a contracting map.
\end{theorem}

\begin{proof}
By Proposition~\ref{pr:Multinucleus}, there exists $n$ such that,
for any adjacency set $A$, the set $(\nuke\cdot A)|_v$ is an
adjacency set for all words $v\in\xs$ of length at least $n$. Fix
such a number $n$ and define then, for $g\in G$ and $v\in\alb^n$,
the point $\wt I_v(g)$ as the barycenter of the simplex
$(\nuke\cdot g)|_v$. For every $h, g\in G$ and $v\in\alb^n$ we
have $(\nuke\cdot h\cdot g)|_v=(\nuke\cdot h)|_{g(v)}\cdot g|_v$,
hence
\[\wt I_v(h\cdot g)=\wt I_{g(v)}(h)\cdot g|_v,\]
i.e., condition~\eqref{eq:phix} is satisfied.

Let $v\in\alb^n$ be an arbitrary word, and let $A\subset G$ be an
adjacency set, i.e., a simplex of $\Gamma$. For every $g\in A$ we
have
\[g|_v\in(\nuke\cdot g)|_v\subset (\nuke\cdot A)|_v,\]
hence all the simplices $(\nuke\cdot g)|_v$ for $g\in A$ belong to
the simplex $\Delta=(\nuke\cdot A)|_v$. Consequently, we can
linearly extend inside $\Delta$ the map $\wt I_v$ from the set of
vertices of the simplex $A$ to the whole geometric realization of
$A$. These extensions agree with each other,
satisfy~\eqref{eq:phix}, and hence define a $G$-equivariant
continuous map $\wt I:\Gamma_n\arr\Gamma$. The points
$I_v(g)=g|_v$ also belong to the simplex $\Delta$, hence the
convex combination $(1-t)I^{(n)}+t\wt I$ inside $\Delta$ is a
$G$-equivariant homotopy from $I^{(n)}$ to $\wt I$.

Fix some $g_0\in A$. For every $g\in A$ we have
$g_0g^{-1}\in\nuke$, hence the simplices $(\nuke\cdot g)|_v$
contain $g_0|_v$ for all $g\in A$. Consequently, their barycenters
$\wt I_v(g)$ are contained in the image of $\Delta$ under the
homothety with center in $I_v(g_0)$ and coefficient
$\dim\Delta/(\dim\Delta+1)$. The maps $\wt I_v$ for
$v\in\bim^{\otimes n}$ are affine on the simplices of $\Gamma$,
hence we get exponential decreasing of the diameters of the
simplices under compositions of the maps $\wt I_v$, for
$v\in\alb^n$, which in turn implies that $\wt I$ is contracting.
\end{proof}

Even though we have proved that only some iteration $(J_0(G),
J_n(G), p^n, \iota^n)$ of the topological nucleus $(J_0(G),
J_1(G), p, \iota)$ is homotopic to a contracting automaton, one
can use this fact to approximate not only the limit space $\lims$,
but also the limit dynamical system $\si:\lims\arr\lims$ (and not
just its $n$th iteration).

Let us introduce some notation. Let $\wt I:\Gamma_n\arr\Gamma$ be
a contracting map equivariantly homotopic to the map
$I^{(n)}:\Gamma_n\arr\Gamma$, as in
Theorem~\ref{th:contractingautomaton}. Denote by
$\wt\iota_k:J_{n(k+1)}(G)\arr J_{nk}(G)$ the map induced by
\[\xi\otimes v_1\otimes v_2\otimes\cdots\otimes
v_k\mapsto\wt I(\xi\otimes v_1)\otimes v_2\otimes\cdots\otimes
v_k,\] for $v_i\in\bim^{\otimes n}$; by
${\wt\iota_k'}:J_{n(k+1)+1}(G)\arr J_{nk+1}(G)$ the map induced by
\[\xi\otimes v_1\otimes v_2\otimes\cdots\otimes v_k\otimes x\mapsto
\wt I(\xi\otimes v_1)\otimes v_2\otimes\cdots\otimes v_k\otimes
x,\] for $v_i\in\bim^{\otimes n}$ and $x\in\bim$; by
$\iota_{nk}:J_{nk+1}(G)\arr J_{nk}(G)$ and $p_{nk}:J_{nk+1}(G)\arr
J_{nk}(G)$, as before, the maps induced by
\[\xi\otimes x\otimes v_1\otimes \cdots\otimes
v_k\mapsto I(\xi\otimes x)\otimes v_1\otimes\cdots\otimes v_k,\]
and
\[\xi\otimes v_1\otimes \cdots\otimes v_k\otimes x\mapsto
\xi\otimes v_1\otimes\cdots\otimes v_k,\] respectively, for
$x\in\bim$ and $v_i\in\bim^{\otimes n}$.

\begin{corollary}
We have two infinite commutative diagrams
\[\begin{array}{cccccc}
\ldots & J_{2n+1}(G) & \stackrel{\wt\iota_1'}{\arr} &
J_{n+1}(G) & \stackrel{\wt\iota_0'}{\arr} & J_1(G)\\
& \mapdown{} & & \mapdown{} & & \mapdown{}\\
\ldots & J_{2n}(G) & \stackrel{\wt\iota_1}{\arr} & J_n(G) &
\stackrel{\wt\iota_0}{\arr} & J_0(G)
\end{array}\]
where in one diagram the vertical arrows are
$\iota_{nk}:J_{nk+1}(G)\arr J_{nk}(G)$, and in the other they are
$p_{nk}:J_{nk+1}(G)\arr J_{nk}(G)$.

The inverse limit of the spaces $J_{nk}$ with respect to the maps
$\wt\iota_k$ and the inverse limit of the spaces $J_{nk+1}$ with
respect to the maps ${\wt\iota_k}'$ are homeomorphic to the limit
space $\lims$. If we identify these limits with each other by the
limit of the maps $\iota_{nk}$, then the limit of the maps
$p_{nk}$ is a dynamical system topologically conjugate to
$\si:\lims\arr\lims$.
\end{corollary}

\begin{proof}
The homeomorphism of the inverse limit of the maps $\wt\iota_k$
with $\lims$ constructed in the proof of
Theorem~\ref{th:approximationlimg} maps the point of the inverse
limit represented by a sequence \[(\xi_1\otimes
v_1,\quad\xi_2\otimes v_2\otimes v_1,\quad\xi_3\otimes v_3\otimes
v_2\otimes v_1, \ldots),\] where $\xi_{k+1}=\wt I(\xi_k\otimes
v_k)$, to the point of $\lims$ represented by $\cdots\otimes
v_3\otimes v_2\otimes v_1\in\limg$. Taking tensor product of the
spaces $\Gamma_{nk}$ with $\bim$, and using the identification
$\limg\cong\limg\otimes\bim$, we get a natural homeomorphism of
the limit of the spaces $\Gamma_{nk+1}$ with the space $\limg$,
mapping
\[(\xi_1\otimes v_1\otimes x,\quad\xi_2\otimes
v_2\otimes v_1\otimes x,\quad\xi_3\otimes v_3\otimes v_2\otimes
v_1\otimes x, \ldots)\] to $\cdots\otimes v_3\otimes v_2\otimes
v_1\otimes x$. It is easily checked now that if we make these
identifications of the inverse limits with $\lims$, then the limit
of the maps $\iota_{nk}$ will be identical on $\lims$, and the
limit of the maps $p_{nk}$ will be the shift $\si:\lims\arr\lims$.
\end{proof}

\section{Examples of contracting topological automata}

\subsection{A self-covering of the torus}
Consider the self-similar group $G$ generated by
\[u=\sigma(1, uv),\quad v=\sigma(u^{-1}, v).\]
It is checked directly that $u$ and $v$ commute and that they have
infinite order. Consequently, this group is isomorphic to the free
abelian group $\Z^2$, and one can apply the general theory
(see~\cite[Section~2.9]{nek:book}) to check that it is
contracting, and to find the limit dynamical system. It is
conjugate to the self-covering of the torus $\C/\Z[i]$ induced by
multiplication by $(1-i)$ on $\C$. Nevertheless, let us apply in
this simple setting the cut-and-paste procedure described in
Proposition~\ref{pr:cutandpaste} and find a simplicial
approximation of the limit dynamical system.

The nucleus of the group generated by $u$ and $v$ is the set
\[\nuke=\{1, u, v, u^{-1}, v^{-1}, uv, u^{-1}v^{-1}\},\]
see~\cite{nek:filling}. Let us find the adjacency sets
$A\subset\nuke$ containing the identity. The set of maximal
simplices of the Rips complex $\Gamma(G, \nuke)$ containing the
identity is
\[\mathcal{A}=\{\{1, u, uv\}, \{1, v, uv\}, \{1, u^{-1}, u^{-1}v^{-1}\},
\{1, v^{-1}, u^{-1}v^{-1}\}, \{1, u, v^{-1}\}, \{1, u^{-1},
v\}\}.\] We have
\begin{alignat*}{2}
\{1, u, uv\}|_0 &= \{1, v\}, &\quad \{1, u, uv\}_1 &= \{1, uv,
v\},\\
\{1, v, uv\}|_0 &= \{1, u^{-1}, v\}, &\quad \{1, v, uv\}|_1 &=
\{1, v\}\\
\{1, u^{-1}, u^{-1}v^{-1}\}|_0 &= \{1, v^{-1}, u^{-1}v^{-1}\},
&\quad \{1, u^{-1}, u^{-1}v^{-1}\}|_1 &= \{1, v^{-1}\},\\
\{1, v^{-1}, u^{-1}v^{-1}\}|_0 &= \{1, v^{-1}\}, &\quad \{1,
v^{-1}, u^{-1}v^{-1}\}|_1 &= \{1, u, v^{-1}\},\\
\{1, u, v^{-1}\}|_0 &= \{1, v^{-1}\}, &\quad \{1, u, v^{-1}\}|_1
&= \{1, u, uv\},\\
\{1, u^{-1}, v\}|_0 &= \{1, u^{-1}, u^{-1}v^{-1}\}, &\quad \{1,
u^{-1}, v\}|_1 &= \{1, v\}.
\end{alignat*}

\begin{figure}
\centering
\includegraphics{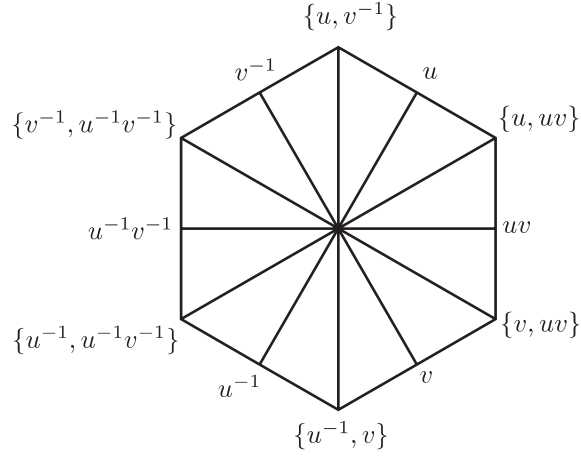}
\label{fig:hexagon}\caption{Complex $T_0$}
\end{figure}

It follows from Corollary~\ref{cor:adjacencygraph} that all the
elements of $\mathcal{A}$ are adjacency sets. The geometric
realization $T_0$ of the poset of adjacency sets containing the
identity is shown on Figure~\ref{fig:hexagon}. (The identity
element is not listed in the labels, e.g., a label $g$ denotes the
vertex corresponding to $\{1, g\}$.)

For every $g\in\nuke\setminus\{1\}$ the set $K_{g, 0}$ is the side
of the hexagon containing the vertex corresponding to $\{1, g\}$.
The transformation $\kappa_{g, 0}$ identifies the opposite sides
$K_{g, 0}$ and $K_{g^{-1}, 0}$ of the hexagon. We conclude that
the complex $J_0(G)$ is a two-dimensional torus.

\begin{figure}
\centering
\includegraphics{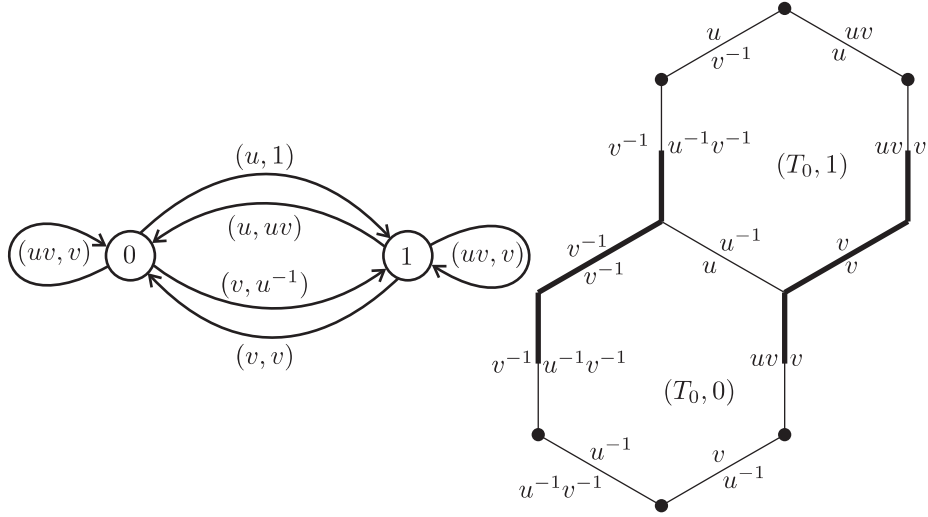}
\label{fig:hexagon2} \caption{The dual Moore diagram of the
nucleus and $T_1$}
\end{figure}

The dual Moore diagram of the nucleus is shown on the left-hand
part of Figure~\ref{fig:hexagon2}. On the right-hand side of the
figure the complex $T_1$ is shown. It follows from the dual Moore
diagram and Proposition~\ref{pr:cutandpaste} that $T_1$ is
obtained by gluing two copies of $T_0$ along two edges (containing
$\{1, u\}$ in the copy $(T_0, 0)$ and containing $\{1, u^{-1}\}$
in $(T_0, 1)$).

The labels inside the hexagons on Figure~\ref{fig:hexagon2}
describe the covering map $p:J_1(G)\arr J_0(G)$, i.e., they repeat
the labels of $T_0$ in its copies $(T_0, x)$. The labels outside
describe the map $\iota:J_1(G)\arr J_0(G)$. The highlighted
vertices are mapped to vertices of the hexagon $T_0$. The labels
show to which sides of the hexagon $T_0$ the corresponding edges
of $T_1$ are mapped (in particular, a letter $g\in\nuke$ labels
the edges of the domains $K_{g, 1}$ of the identifications
$\kappa_{g, 1}$, described in Proposition~\ref{pr:cutandpaste}).

\begin{figure}
\centering
\includegraphics{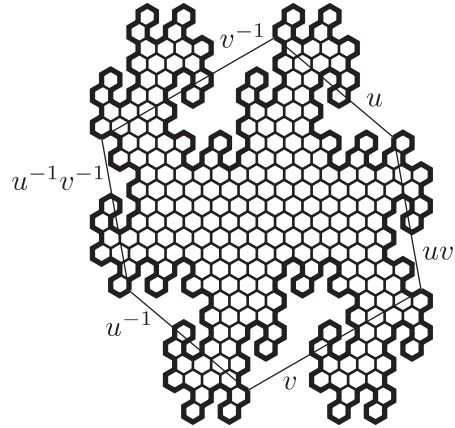}
\label{fig:dragon} \caption{The complex $T_6$}
\end{figure}

The simplicial map $\iota$ maps the two highlighted portions of
$T_1$ to single vertices ($\{1, v\}$ and $\{1, v^{-1}\}$,
respectively): it maps the top half of the top hexagon and the
bottom half of the bottom hexagon to the top and the bottom halves
of the hexagon $T_0$, respectively; the remaining part of $T_1$ is
mapped to the horizontal axis of symmetry of $T_0$ (passing
through $\{1, uv\}$ and $\{1, u^{-1}v^{-1}\}$. See
Figure~\ref{fig:dragon}, where the complex $T_6$ is shown, which
was obtained by application of Proposition~\ref{pr:cutandpaste}.
The hexagon $T_0$ is superimposed with $T_6$ in such a way that
the vertices of the hexagon $T_0$ coincide with their preimages
under $\iota^6$.

It is easy to see that the map $\iota:J_1(G)\arr J_0(G)$ is
homotopic to a homeomorphism $\wt\iota$. When we replace $\iota$
by $\wt\iota$, we will transform the topological nucleus of $G$
into a subdivision rule $(J_0(G), J_1(G), p, \wt\iota)$ defining a
self-covering $p:J_1(G)\arr J_0(G)$ of a torus (where $J_1(G)$ and
$J_0(G)$ are identified with each other by $\wt\iota$). For
instance, by computing the action of the self-covering $p$ on the
homology, we conclude that it is homotopic to the self-covering of
$\C/\Z[i]$ induced by $z\mapsto (1-i)z$.

\subsection{A Forn\ae ss-Sibony example}
The following rational transformation of $\C^2$ was studied
in~\cite{fornsibon:crfin}.
\[f\left(z, p\right)=
\left(\left(1-\frac{2z}{p}\right)^2,\quad
\left(1-\frac{2}{p}\right)^2\right).\]

The map $f$ can be extended to an endomorphism of $\CP^2$. The
post-critical set of $f$ is then the union of the lines $z=1, z=0,
p=1, p=0$, $p=z$ and the line at infinity.

The iterated monodromy group of $f$ was computed
in~\cite{nek:dendrites}. It is easier to describe the iterated
monodromy group of the quotient of the map $f$ by the complex
conjugation $(z, p)\mapsto (\overline z, \overline p)$, which will
be an index two extension of $\img{f}$. It is the group $G$
generated by the transformations
\begin{alignat*}{2}
\alpha &= \sigma, &\qquad  a &=\pi,\\
\beta  &=\left(\alpha, \gamma, \alpha, \gamma\right), &\qquad  b
&=\left(a\alpha, a\alpha, c, c\right),\\
\gamma &=\left(\beta, 1, 1, \beta\right), &\qquad  c
&=\left(b\beta, b\beta, b, b\right),
\end{alignat*}
where $\sigma=(12)(34)$ and $\pi=(13)(24)$. The iterated monodromy
group of $f$ is generated then by $\alpha, \beta, \gamma,
S=ac\gamma$ and $T=cb$.

Direct computation shows that the nucleus of $G$ is a union of the
following six finite groups
\begin{eqnarray*}
G_A= & \langle\beta, \gamma, b, c\rangle & \cong D_8\rtimes D_4,\\
G_B= & \langle\alpha, \gamma, a, c\rangle & \cong D_4\rtimes
D_2,\\
G_C= & \langle\alpha, \beta, a, b\rangle & \cong D_8\rtimes D_4,
\end{eqnarray*}
and
\begin{eqnarray*}
G_{A_1}= &\langle\alpha, b, c\rangle& \cong
C_2\times D_4,\\
G_{B_1}= &\langle\beta, a, c\gamma\rangle
 & \cong C_2\times D_2,\\
G_{C_1}= &\langle\gamma, a\alpha, b\beta\rangle & \cong C_2\times
D_4,
\end{eqnarray*}
where $D_n$ denotes the dihedral group of order $2n$ and $C_n$ is
a cyclic group of order $n$. Note that the group of inner
automorphisms of $D_{2n}$ is isomorphic to $D_n$, which defines
the corresponding semidirect products above.

Inspection of the Moore diagram of the nucleus shows that these
six subgroups $G_*$ are precisely the maximal adjacency sets
containing the identity element.

Consider the poset $\mathfrak{G}$ of the subgroups $G_*$ and their
all possible intersections (pairwise and triple are enough, since
all the rest are trivial). One can show that for every
$H\in\mathfrak{G}$ and every $x\in\{1, 2, 3, 4\}$, the set $H|_x$
is also an element of $\mathfrak{G}$. It follows that the set
$\overline\Gamma$ of cosets $H\cdot g$ for $H\in\mathfrak{G}$ and
$g\in G$ is a $G$-invariant subcomplex of the barycentric
subdivision $\Gamma'$ of the tiling nerve $\Gamma$ of $G$, and
that $\overline\Gamma$ is invariant under the maps
$I_x:\Gamma'\arr\Gamma'$.

It follows that restricting the equivariant map
$I:\Gamma\otimes\bim\arr\Gamma$ onto $\overline\Gamma\otimes\bim$
we get a combinatorial model $(\overline J_0(G), \overline J_1(G),
p, \iota)$ of $f$.

The complex $\overline J_0(G)$ is the geometric realization of the
poset $\mathfrak{G}$. It is a union of three tetrahedra with a
common face.

The recursive definition of the complexes $\overline J_n(G)$,
approximating the limit space of $G$ is a simple pasting rule,
which has a nice interpretation in the spirit of Hubbard trees.
The Julia set of $f$ is approximated then by two copies of
$\overline J_n(G)$ glued together in a natural way. See for more
details the paper~\cite{nek:dendrites}.

\subsection{Post-critically finite rational functions}

Let $f:\CS\arr\CS$ be a post-critically finite complex rational
function. Let $\M$ be the Thurston orbifold of $f$ and let
$\F=(\M, \M_1, f, \iota)$ be the associated topological automaton.
The underlying space of the orbifold $\M$ is a punctured sphere if
$f$ has a super-attracting cycle (i.e., a cycle containing a
critical point) and is the whole sphere, if every critical point
is strictly pre-periodic. In all these cases $f$ is expanding with
respect to the Poincar\'e metric on $\M$. There exists a compact
subset $\M'\subset\M$ such that it contains all singular points of
$\M$, and $f^{-1}(\M')\subset\M'$ (one can take $\M'$ to be the
set bounded by appropriate level curves of the Green function of
the Julia set of $f$). Restricting the Poincar\'e metric onto
$\M'$ we get a contracting topological automaton $\F'=(\M', \M_1',
f, \iota)$.

\begin{proposition}
If $f$ has a super-attracting cycle (in particular, when it is a
polynomial), then there exists $n$ such that $f^{\circ n}$ is
combinatorially equivalent to a contracting automaton $\F=(J, J_1,
f, \iota)$, where $J$ is a graph of cyclic groups.
\end{proposition}

\begin{proof}
If $f$ has a super-attracting cycle, then the orbifold $\M'$ can
be retracted to a graph of groups, and we can use
Proposition~\ref{pr:simplification}.
\end{proof}

In many cases, choosing a nice retract of the Thurston orbifold and
choosing a correct metric on the retract, one can find a
contracting combinatorial model of $f$, and not just of $f^{\circ
n}$ for some $n$. Such constructions are classical for
post-critically finite polynomials.

\subsection{Hubbard trees of strictly pre-periodic polynomials}

Let $f$ be a post-critically finite polynomial such that every
finite critical point of $f$ is strictly pre-periodic (i.e., has
finite forward orbit, but does not belong to a cycle). Then the
polynomial has no attracting cycles in $\C$, therefore its Julia
set $J_f$ is a dendrite, i.e., every two points $x, y\in J_f$ can
be connected by a unique arc. The set $P_f\setminus\{\infty\}$ of
finite post-critical points of $f$ is a subset of $J_f$. The
\emph{Hubbard tree} of $f$ is the convex hull of the set of finite
post-critical points in $J_f$. Here the convex hull of a set
$A\subset J_f$ is the union of the arcs connecting all pairs of
points of $A$. The Hubbard tree is a natural choice for a retract
of the Thurston orbifold.

The Hubbard tree $H_f$ is invariant, i.e., $f(H_f)=H_f$. For every
point $x\in J_f$ there exists a unique point $y\in H_f$ such that
the arc connecting $x$ with $y$ has no common points with $H_f$
except for $y$. The point $y$ is called \emph{projection} of $x$
onto $H_f$. In particular, projection of a point $x\in H_f$ onto
$H_f$ is the point $x$ itself. It is not hard to show that the
projection map $\iota:J_f\arr H_f$ is continuous.

Denote by $\M$ the orbispace with the underlying space $H_f$, with
the orbispace structure obtained by restricting the Thurston
orbispace of $f$ onto $H_f$ (see the definition of the Thurston
orbispace in Subsection~\ref{sss:thurstonmaps}).
Let $\M_1$ be the orbispace with the underlying space
$f^{-1}(H_f)$ such that $f:\M_1\arr\M$
is a $\deg(f)$-fold covering of orbispaces. Then the projection
map $\iota:f^{-1}(H_f)\arr H_f$ is a morphism of the orbispaces.

The iterate $\M_n$ of the automaton $(\M, \M_1, f, \iota)$ is
homeomorphic to the convex hull of the set $f^{-n}(P_f)$ in the
Julia set $J_f$.

The obtained topological automaton is contracting with respect to
an appropriate metric on $H_f$. It is combinatorially equivalent
to the polynomial $f$ (i.e., to the corresponding partial
self-coverings) by Proposition~\ref{pr:combequiv}.

Consequently, the Hubbard tree is a model of the dynamical system
$(J_f, f)$, by Theorem~\ref{th:projlim}. Hubbard trees are used
extensively in symbolic dynamics of polynomial iterations,
see~\cite{DH:orsayI,DH:orsayII,henkdierk}.

\begin{figure}
\centering  \includegraphics{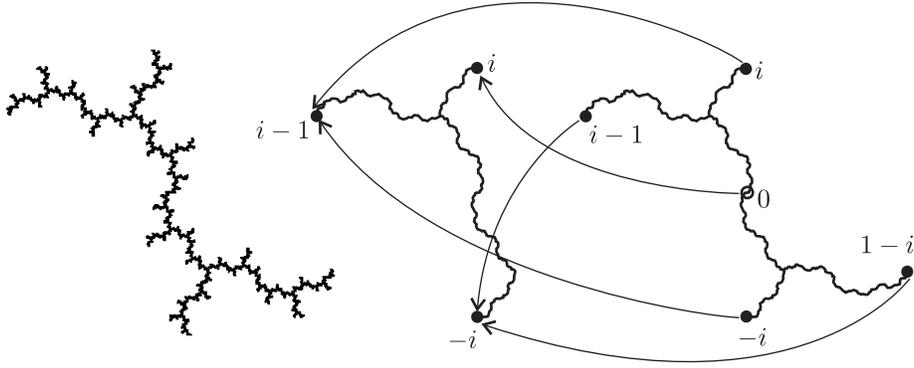}\\
  \caption{Hubbard tree of $z^2+i$}\label{fig:i}
\end{figure}

As an example consider the polynomial $z^2+i$. The orbit of its
critical value $i$ is $i\mapsto -1+i\mapsto -i\mapsto -1+i$.
Figure~\ref{fig:i} shows on its left-hand side the Julia set of
$z^2+i$. On the right-hand side of the figure the Hubbard trees
$\M$ and $\M_1$ are shown. Black dots mark the singular points of
the corresponding orbispaces. Isotropy group of each of the
singular points is of order 2. The point 0 is critical (but
non-singular). The morphism $\iota:\M_1\arr\M$ maps the whole
branch of $\M_1$ containing $1-i$ to the point at which this
branch is connected to $\M$, and ``kills'' the isotropy group. It
acts identically on the rest of the Hubbard tree.

\subsection{Hubbard graphs of polynomials}

Hubbard trees can be also defined for arbitrary post-critically
finite polynomials, see~\cite{DH:orsayI,DH:orsayII}. A more
appropriate construction in our setting is a modification of the
classical construction of Hubbard trees, obtained by replacing some
of its vertices by circles. Instead of formulating a general
construction, we will just describe two examples.

\subsubsection{Basilica}
\begin{figure}
\centering
  \includegraphics{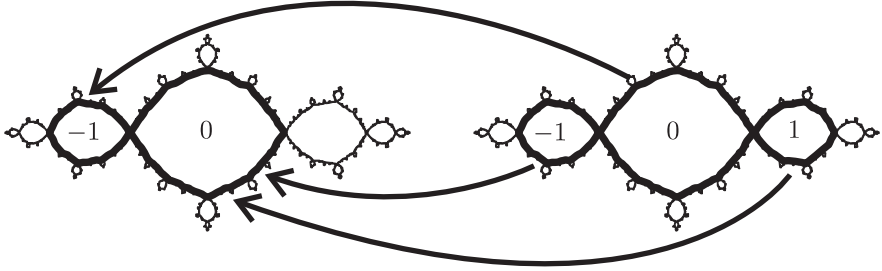}\\
  \caption{Basilica}\label{fig:basilicaj}
\end{figure}

Consider the polynomial $f(z)=z^2-1$. Its post-critical set is
$\{0, -1, \infty\}$. The Julia set of $z^2-1$, called
\emph{Basilica}, is shown on Figure~\ref{fig:basilicaj}. Let $\M$
be the union of the boundaries of the Fatou components containing
the finite post-critical points $0$ and $-1$ (it is highlighted on
the left-hand side part of the figure). The set $\M$ is
homotopically equivalent to $\widehat\C\setminus\{0, -1, \infty\}$
and is forward invariant. Let $\M_1=f^{-1}(\M)$. It is the union
of the boundaries of the Fatou components of $0, 1$ and $-1$
(highlighted on the right-hand side part of
Figure~\ref{fig:basilicaj}). The arrows on
Figure~\ref{fig:basilicaj} show the action of $f$.

It is easy to see now that the topological automaton
\[\F=(\C\setminus\{0, -1\}, \C\setminus\{0, 1, -1\}, f, id)\]
is homotopically equivalent to the automaton $\F_1=(\M, \M_1, f,
\iota)$, where $\iota$ is identical on $\M\subset\M_1$ and maps
the boundary of the Fatou component of $1$ to its common point
with $\M$.

The space $\M_n$ obtained by iteration of the constructed
automaton is homeomorphic to the union of the boundaries of the
Fatou components containing the points of the set $f^{-n}(\{0, -1\})$.

\begin{figure}
\centering\includegraphics{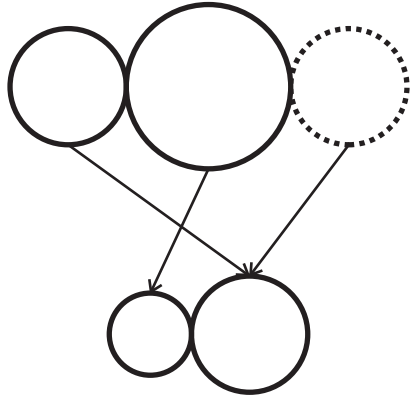}\\
\caption{Affine model of Basilica}\label{fig:basilicam}
\end{figure}

The fact that the automaton $\F_1$ is contracting can be shown
directly by introduction of a natural length structure on $\M$ and
$\M_1$, and considering an abstract affine model of $\F_1$. The
space $\M$ will be a one point union of a circle of length 1 and a
circle of length $\sqrt{2}$. Let $\M_1$ be the double locally
isometric covering of $\M$ such that the circle of length 1 is
doubly covered by a circle of length 2, and the circle of length
$\sqrt{2}$ is covered by two isometric circles. See the covering
on Figure~\ref{fig:basilicam}.

Let $\iota:\M_1\arr\M$ be the continuous map contracting one of
the two $f$-preimages of the circle of length $\sqrt{2}$ (shown by
a dashed line on Figure~\ref{fig:basilicam}) to its common point
with the circle of length $2$, dividing by $\sqrt{2}$ all the
distances in the other two circles of $\M_1$ and then mapping them
isometrically onto $\M$. The obtained automaton is topologically
conjugate to $\F_1$.

\subsubsection{Airplane}

\begin{figure}
\centering\includegraphics{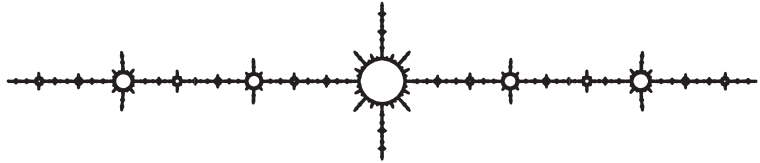}\\
\caption{Airplane}\label{fig:airplane}
\end{figure}

Figure~\ref{fig:airplane} shows the Julia set the ``Airplane''
polynomial $z^2+c$ for $c\approx -1.7549\ldots$, which is
determined by the condition that it has real coefficients, and the
critical point 0 belongs to a cycle of length three. We can use
again the boundaries of the Fatou components of the post-critical
points to construct a simple contracting topological automaton
combinatorially equivalent to the polynomial. The main difference
with the case of the polynomial $z^2-1$ is that these boundaries
are disjoint.

\begin{figure}
\centering\includegraphics{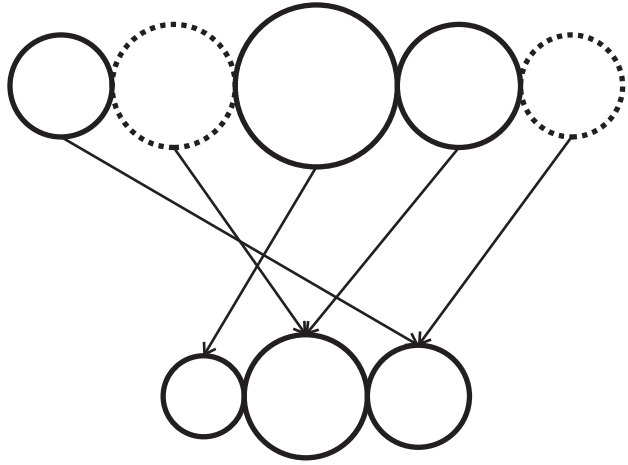}\\
\caption{Model of the Airplane} \label{fig:airplane2}
\end{figure}

Hence, one has to attach the circles corresponding to the boundaries to
each other imitating their relative arrangement in the Julia set.
The corresponding topological automaton is shown on
Figure~\ref{fig:airplane2}.

The circles of $\M$ are of lengths $1, \sqrt[3]{4}, \sqrt[3]{2}$;
the space $\M_1$ covers the circle of the length $1$ twice by a
circle of length $2$, and covers the remaining two circles
isometrically by pairs of circles. The dashed lines show the
circles, which are collapsed by the map $\iota$. It divides the
lengths of the other circles by $\sqrt[3]{2}$.

\subsection{Correspondences on moduli spaces}
Let $f:S^2\arr S^2$ be a Thurston map (i.e., an orientation
preserving post-critically finite branched self-covering of the sphere) of
degree $d$. Let $P_f$ be the post-critical set of $f$.

We present here a short summary of the Teichm\"uller theory of
Thurston maps. For more details and for relation of these concepts
with a theorem of Thurston, see~\cite{DH:Thurston}.

The \emph{Teichm\"uller space} $\mathcal{T}_{P_f}$ modelled on
$(S^2, P_f)$ is the space of homeomorphisms $\tau:S^2\arr\CS$
(seen as complex structures on $S^2$), where two complex
structures $\tau_1, \tau_2:S^2\arr\CS$ are identified if there
exists a M\"obius transformation $\phi:\CS\arr\CS$ such that
$\phi\circ\tau_1$ is isotopic to $\tau_2$ relative to $P_f$ (and
is equal to $\tau_2$ on $P_f$).

For every complex structure $\tau\in\mathcal{T}_{P_f}$ there
exists a unique complex structure $\tau'\in\mathcal{T}_{P_f}$,
such that the map $f_\tau=\tau\circ f\circ(\tau')^{-1}$ closing
the commutative diagram
\begin{equation}\label{eq:sigmaback}\begin{array}{ccc}S^2 & \stackrel{f}{\arr} & S^2\\
\mapdown{\tau'} & & \mapdown{\tau}\\
\CS & \stackrel{f_\tau}{\arr} & \CS\end{array}\end{equation} is a
rational function. Let us denote $\tau'=\sigma_f(\tau)$.

The \emph{moduli space} $\M=\M_{P_f}$ of $(S^2, P_f)$ is the space
of injective maps $P_f\arr\CS$ modulo compositions with M\"obius
transformations. It is known that $\mathcal{T}_{P_f}$ is the
universal covering of $\M_{P_f}$, where the covering map is
$\tau\mapsto\tau|_{P_f}$.

The fundamental group of $\M_{P_f}$ can be identified with the
(pure) mapping class group $G$ of $(S^2, P_f)$, so that the action
of the fundamental group on the universal covering
$\mathcal{T}_{P_f}$ coincides with the action of the mapping class
group on $\mathcal{T}_{P_f}$ by compositions with the
homeomorphisms:
\[g:\tau\mapsto\tau\circ g\]
for $g\in G$ and $\tau\in\mathcal{T}_{P_f}$.

Every homeomorphism $g\in G$ can be lifted by the Thurston map
$f:S^2\arr S^2$ to a homeomorphism of $S^2$, which we will denote
by $f^*g$. Let $G_1$ be the subgroup of the elements $g\in G$ such
that $f^*g$ fixes the set $P_f\subset f^{-1}(P_f)$ pointwise,
i.e., is an element of $G$. It is easy to see that $G_1$ is a
subgroup of finite index in $G$. The map $g\mapsto f^*g$ is a
homomorphism from $G_1$ to $G$, i.e., it is a virtual endomorphism
of the group $G$.

Let $\M_1$ be the quotient of $\mathcal{T}_{P_f}$ by the action of
$G_1$. Then the identity map on $\mathcal{T}_{P_f}$ induces a
finite covering $F:\M_1\arr\M$.

On the other hand, since $f^{-1}(P_f)\supseteq P_f$, we get a well
defined continuous map
\[\iota:\M_1\arr\M:\tau\mapsto\sigma_f(\tau)|_{P_f}.\]

We will call the obtained topological automaton $(\M, \M_1, F,
\iota)$ the \emph{moduli space correspondence} associated with the
Thurston map $f$. It is easy to check that the virtual
endomorphism $g\mapsto f^*g$ is associated with the topological
automaton $(\M, \M_1, F, \iota)$.

As before, we interpret the topological automaton as the
correspondence $\iota(x)\mapsto f(x)$, which in this case is the
projection of the correspondence $\sigma_f(\tau)\mapsto\tau$ onto
the moduli space.

In some cases (but not in general) $\iota$ is one-to-one and we
get hence partial self-covering of $\M$.

As an example consider a quadratic polynomial $f(z)=z^2+c$ such
that the critical point $0$ belongs to a cycle of length $n$, and
look at it just as at a Thurston map $f:S^2\arr S^2$. Let $\infty,
0, c=z_0, \ldots, z_{n-2}$ be the post-critical set $P_f$, where
$z_k=f(z_{k-1})$ and $0=f(z_{n-2})$. Let $\tau:S^2\arr\CS$ be an
arbitrary point of the Teichm\"uller space $\mathcal{T}_{P_f}$.
The corresponding point of the moduli space $\M$ is determined by
the values of $\tau(\infty)$ and $\tau(z_k)$ for $k=0, \ldots,
n-2$.

Applying an appropriate M\"obius transformation, we may assume
that $\tau(\infty)=\infty$, $\tau(0)=0$ and $\tau(z_0)=1$. It
follows that the corresponding point of the moduli space is
determined by the tuple
\[(\tau(z_1), \tau(z_2), \ldots, \tau(z_{n-2}))=(p_1, p_2, \ldots,
p_{n-2})\in\C^{n-2}.\] In this way we identify the moduli space
$\M$ with the set
\[\{(p_1, p_2, \ldots, p_{n-2})\in\C^{n-2}\;:\;p_i\ne 0, p_i\ne 1,
p_i\ne p_j\text{\ for\ }i\ne j\}.\]

Let $(p_1', p_2', \ldots, p_{n-2}')\in\M$ be the point of the
moduli space corresponding to $\sigma_f(\tau)$. Then it follows
from the commutative diagram~\eqref{eq:sigmaback}, that $f_\tau$
is a quadratic polynomial with critical point $0$ such that
\[f_\tau(0)=1, f_\tau(1)=p_1', f_\tau(p_1)=p_2', \ldots,
f_\tau(p_{n-3})=p_{n-2}', f_\tau(p_{n-2})=0.\]

The first equality and the fact that $0$ is a critical point imply that
$f_\tau(z)=1+az^2$ for some $a\in\C$. It follows from the equality
$f_\tau(p_{n-2})=0$ that $a=-\frac{1}{p_{n-2}^2}$. The remaining
equalities imply
\[(p_1', p_2', \ldots, p_{n-2}')=\left(1-\frac{1}{p_{n-2}^2},
1-\frac{p_1^2}{p_{n-2}^2}, \ldots,
1-\frac{p_{n-3}^2}{p_{n-2}^2}\right).\]

It follows that the correspondence $\sigma_f(\tau)\mapsto\tau$ is
projected onto the moduli space to the rational function $(p_1,
\ldots, p_{n-2})\mapsto(p_1', \ldots, p_{n-2}')$, so that the
moduli space correspondence is a partial self-covering.

Note that in this case the rational map can be extended to an
endomorphism of $\CP^{n-2}$ given in homogeneous coordinates by
\[[p_1:p_2:\cdots:p_{n-2}:p_{n-1}]\mapsto
[p_{n-2}^2-p_{n-1}^2:p_{n-2}^2-p_1^2:\cdots:p_{n-2}^2-p_{n-3}^2:p_{n-2}^2].\]
The post-critical set of this endomorphism is the union of the
lines $p_i=0, p_i=1, p_i=p_j$ for all $i=1, \ldots, n-1$ and $i\ne
j=1, \ldots, n-1$.

For more on this and similar post-critically finite endomorphisms
of complex projective spaces, see~\cite{koch:french}.

Let us describe a combinatorial model of the moduli space
correspondence for the case when $f$ is a topological polynomial.
Here a Thurston map $f:S^2\arr S^2$ is called a topological
polynomial if there exists a point $x\in S^2$ such that
$f^{-1}(x)=\{x\}$. In this case $x$ is called the \emph{point at
infinity}, and $f$ is considered to be a branched self-covering of
the plane $S^2\setminus\{x\}$.

Our combinatorial model will be a topological automaton $(D_n,
D_n', p, \iota)$, where $D_n$ and $D_n'$ are affine polyhedra, the
covering $p$ is a local isometry and $\iota$ is contracting,
provided the topological polynomial $f$ is \emph{hyperbolic},
i.e., every cycle of $f$ in $P_f$ contains a critical point. The
space $D_n$ is the moduli space of a family of planar graphs, thus
its definition is similar to the construction of the classifying
space of the outer automorphism group of the free group, given
in~\cite{cullervogtmann:moduli} by M.~Culler and K.~Vogtmann. The
complex $D_n$ is also closely related to the classifying space of
the braid groups defined in~\cite{brady:newkp1}.

The polyhedron $D_n$ will depend only on the size of the
post-critical set $P_f$. Let $|P_f|=n+1$ so that $f$ has $n$
finite post-critical points. A \emph{cactus diagram of $n$ discs}
is an oriented two-dimensional contractible cellular complex
$\Gamma$ consisting of $n$ discs labeled by numbers from $1$ to
$n$, such that any two disc are either disjoint or have only one
common point on their boundaries. A \emph{planar cactus diagram}
is a cactus diagram together with an isotopy class of an
orientation preserving embedding $\Delta:\Gamma\arr\R^2$ into the
plane. The isotopy class is uniquely determined by the cyclic
orders of the discs adjacent to every given disc of the diagram.

A \emph{metric cactus diagram} is a cactus diagram together with a
metric on the one-skeleton of the diagram, such that perimeter of
a disc labeled by $k$ is equal to a fixed positive number $l_k$. A
\emph{planar metric diagram} is a metric cactus diagram together
with an isotopy class of an orientation preserving embedding into
the plane.

The cells of the polyhedron $D_n$ are in a bijective
correspondence with the planar cactus diagrams of $n$ discs, while
the points of $D_n$ are in a bijective correspondence with metric
planar cactus diagrams. Points of a given cell are obtained by
specifying the lengths of the edges in the one-skeleton of the
corresponding diagram, so that the perimeters of the discs are
equal to the chosen numbers $l_k$. It follows that dimension of a
cell is equal to the number of the vertices of the corresponding
diagram minus one. When some of the distances go to zero, the
number of vertices of the planar diagram decreases and the
corresponding point of $D_n$ approaches to a cell of lower
dimension.

In particular, the polyhedron $D_n$ has $(n-1)!$ vertices,
corresponding to planar diagrams in which all discs have one
common point (a bouquet of discs). There are no distances to
specify. One-dimensional edges of $D_n$ correspond to diagrams
with two vertices, so that we have to specify one distance. The
maximal number of vertices for a cactus diagram of $n$ discs is
$n-1$, so the polyhedron $D_n$ is $(n-2)$-dimensional.

We use the lengths of the edges of the one-skeleta of the diagrams
as affine coordinates on the corresponding cell. For a given
planar diagram, the set of possible metric realizations (i.e., the
corresponding cell) is a direct product of simplices, due to the
constrains on the perimeters of each of the discs.

\begin{figure}
\centering
\includegraphics{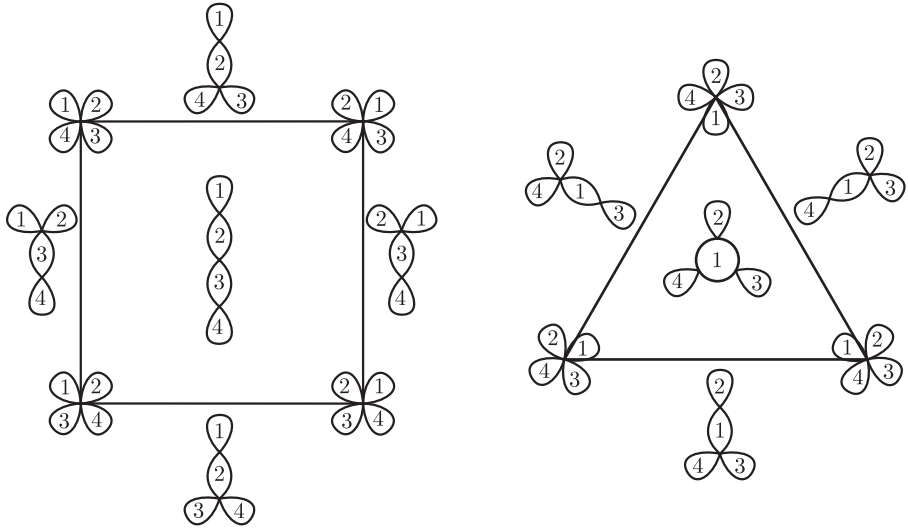}
\caption{Cells of $D_4$}\label{fig:faces}
\end{figure}

\begin{figure}
\centering
\includegraphics{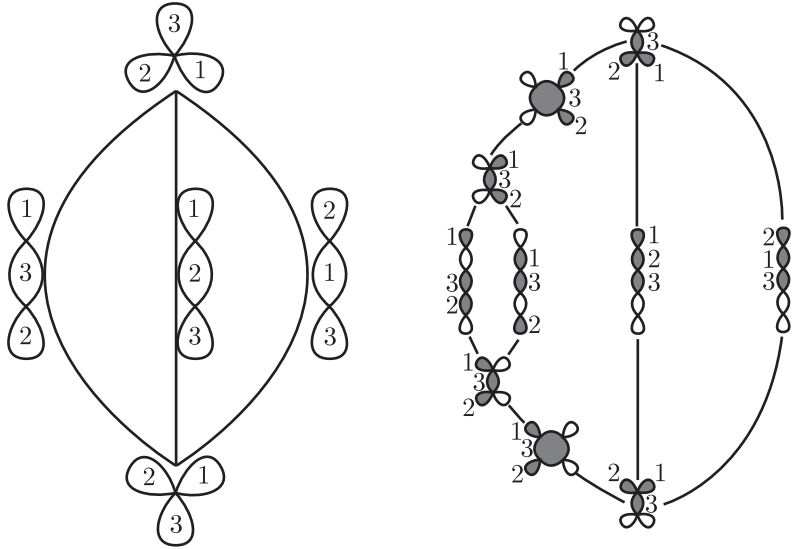}\\
\caption{A moduli space model of $1-1/z^2$}\label{fig:rabbit}
\end{figure}

As an example, see Figure~\ref{fig:faces}, where 2-cells of $D_4$
are described. The complex $D_3$ is shown on the left-hand side of
Figure~\ref{fig:rabbit}.

Let now $f:\R^2\arr\R^2$ be a topological polynomial with $n$
post-critical points $z_1, \ldots, z_n$.

For every metric planar diagram $\Delta\in D_n$ find a
representative $\wt\Delta:\Gamma\arr\R^2$ of the corresponding
isotopy class such that the image of the disc labeled by $k$
contains in its interior the point $z_k$ for every $k=1, \ldots,
n$. Let $f^{-1}(\wt\Delta)$ be the lift of $\wt\Delta$ by $f$. It
is a planar diagram of $|f^{-1}(\{z_1, \ldots, z_n\})|$ discs with
one preimage of a post-critical point in each disc. We introduce a
metric on this planar diagram by lifting it from $\Delta$. We have
$\{z_1, \ldots, z_n\}\subset f^{-1}(\{z_1, \ldots, z_n\})$, so
that all post-critical points belong to interiors of some of the
discs of the diagram $f^{-1}(\wt\Delta)$. For each given $\Delta$
there is only a finite number of possibilities for the isotopy
class of $f^{-1}(\wt\Delta)$ and for the assignments of the
post-critical points to the discs of $f^{-1}(\wt\Delta)$. We get
in this way a finite number of metric planar diagrams
$f^{-1}(\wt\Delta)$ in which some discs are labeled by
post-critical points $z_k$. The space $D_n'$ of such diagrams is
also an affine polyhedron, such that the map
$p:f^{-1}(\wt\Delta)\mapsto\Delta$ is an isometric covering.

In each of the labeled diagrams $f^{-1}(\wt\Delta)$ contract the
non-labeled discs to points and rescale the perimeters of the
remaining labeled discs so that the disc containing the point
$z_k$ has perimeter $l_k$. We will get in this way a metric planar
diagram $\iota(f^{-1}(\wt\Delta))\in D_n$
(we label the disc containing $z_k$ by $k$). We have
defined in this way a piecewise affine map $\iota:D_n'\arr D_n$.

\begin{theorem}
The topological automaton $\F=(D_n, D_n', p, \iota)$ associated
with a post-critically finite polynomial $f$ is combinatorially
equivalent to the moduli correspondence associated with $f$.

If the polynomial $f$ is hyperbolic then the automaton $\F$ is
contracting.
\end{theorem}

\begin{proof}
Let $P=\{z_i\}_{i=1,\ldots, n}$ be, as above, the set of finite
post-critical points of $f$. Consider the space $\wt{D}_n$ of
isotopy classes relative to $P$ of embeddings
$\Delta:\Gamma\arr\R^n$ of metric cactus diagrams of $n$ discs
such that the image of the disc number $k$ contains the point
$z_k$ in its interior.

The pure mapping class group $G$ of $(S^2, P\cup \infty)$ acts
naturally on $\wt{D}_n$. The action is free (in particular, since
the action of the mapping class group on the fundamental group of
$\R^2\setminus P$ by outer automorphism group is faithful). The
quotient of $\wt{D}_n$ by the action is the space $D_n$, hence the
action is co-compact.

For an embedding $\Delta:\Gamma\arr\R^n$, representing a point of
$\wt{D}_n$, consider the lift of $\Delta$ by $f$, contract in the
lift the discs that do not contain points of $P$ in their
interior, and rescale the perimeters of the remaining discs
accordingly to the indices of post-critical points contained in
them. We will get then a point $\Phi(\Delta)$ of $\wt{D}_n$. The
map $\Phi$ obviously satisfies the condition
\[\Phi(\Delta\cdot g)=\Phi(\Delta)\cdot (f^*g),\]
for all elements $g\in G_1$. Here $\Delta\cdot g$ is the image of
$\Delta$ under the action of $g$, and $G_1\le G$ is the subgroup
of elements of $G$ lifted to elements of $G$ by the branched
covering $f$ (see above). It follows that $\Phi$ agrees with the
virtual endomorphism associated with the moduli correspondence
$(\M, \M_1, F, \iota)$ associated with $f$, hence the associated
topological automaton $\mathcal{D}=(\wt{D}_n/G, \wt{D}_n/G_1, P,
\varphi)$, where $P:\wt{D}_n/G_1\arr\wt{D}_n/G$ is the covering
induced by the inclusion $G_1<G$, and $\varphi$ is the map induced
by $\Phi$, is combinatorially equivalent to the moduli
correspondence.

It follows directly from the definitions that the topological
automaton $\mathcal{D}$ is isomorphic to $(D_n, D_n', p, \iota)$.
It is also easy to show that if every cycle of $f$ contains a
critical point, then $\Phi$ is contracting.
\end{proof}

Figure~\ref{fig:rabbit} shows the complexes $D_3$ and $D_3'$ for a
quadratic polynomial $f$ such that its finite critical point
belongs to a cycle of length three. The labels $1, 2, 3$
correspond to the post-critical points $z_1, z_2, z_3$, where
$f(z_k)=z_{k+1}$ and $z_3$ is the critical point. On the right
hand side of Figure~\ref{fig:rabbit} we show the diagrams
$f^{-1}(\wt\Delta)$, where grey cells are the cells containing
post-critical points.

\providecommand{\bysame}{\leavevmode\hbox
to3em{\hrulefill}\thinspace}


\begin{thebibliography}{Nek08b}

\bibitem[Ale29]{alexandroff:simplicialapprox}
P.~Alexandroff, \emph{{\"Uber Gestalt und Lage abgeschlossener
Mengen
  beliebiger Dimension}}, Annals of Mathematics (2) \textbf{30} (1928--1929),
  101--187.

\bibitem[Ale83]{al:free_en}
S.~V. Aleshin, \emph{A free group of finite automata}, Moscow
University
  Mathematics Bulletin \textbf{38} (1983), 10--13.

\bibitem[BGN03]{bgn}
Laurent Bartholdi, Rostislav Grigorchuk, and Volodymyr
Nekrashevych, \emph{From
  fractal groups to fractal sets}, {Fractals in Graz 2001. Analysis -- Dynamics
  -- Geometry -- Stochastics} (Peter Grabner and Wolfgang Woess, eds.),
  {Birkh\"auser Verlag, Basel, Boston, Berlin}, 2003, pp.~25--118.

\bibitem[BH99]{bridhaefl}
Martin~R. Bridson and Andr\'e Haefliger, \emph{Metric spaces of
non-positive
  curvature}, Grundlehren der Mathematischen Wissenschaften, vol. 319,
  Springer, Berlin, 1999.

\bibitem[BL01]{lub:treelatices}
Hyman Bass and Alexander Lubotzky, \emph{Tree lattices}, Progress
in
  Mathematics, vol. 176, {Birkh\"auser Boston Inc.}, {Boston, MA}, 2001, {With
  appendices by Bass, L. Carbone, Lubotzky, G. Rosenberg and J. Tits}.

\bibitem[BP94]{bullett:gallery}
Shaun Bullett and Christopher Penrose, \emph{A gallery of iterated
  correspondences}, Experiment. Math. \textbf{3} (1994), no.~2, 85--105.

\bibitem[Bra70]{brauer:french}
Wilfried Brauer, \emph{Automates topologiques et ensembles
reconnaissables},
  S\'eminaire M. P. Sch\"utzenberger, A. Lentin et M. Nivat, 1969/70:
  Probl\`emes Math\'ematiques de la Th\'eorie des Automates, Secr\'etariat
  math\'ematique, Paris, 1970, pp.~Exp. 18, 24.

\bibitem[Bra01]{brady:newkp1}
Thomas Brady, \emph{A partial order on the symmetric group and new
  {$K(\pi,1)$}'s for the braid groups}, Adv. Math. \textbf{161} (2001), no.~1,
  20--40.

\bibitem[BS02]{henkdierk}
Henk Bruin and Dierk Schleicher, \emph{Symbolic dynamics of
quadratic
  polynomials}, Institut Mittag-Leffler, Report No. 7, 2001/2002.

\bibitem[Bul88]{bullett:quadratic}
Shaun Bullett, \emph{Dynamics of quadratic correspondences},
Nonlinearity
  \textbf{1} (1988), no.~1, 27--50.

\bibitem[Bul91]{bullett:agm}
\bysame, \emph{Dynamics of the arithmetic-geometric mean},
Topology \textbf{30}
  (1991), no.~2, 171--190.

\bibitem[Bul92]{bullet:crfinitemodular}
\bysame, \emph{Critically finite correspondences and subgroups of
the modular
  group}, Proc. London Math. Soc. (3) \textbf{65} (1992), no.~2, 423--448.

\bibitem[CFP01]{cfp:subdivision}
J.~W. Cannon, W.~J. Floyd, and W.~R. Parry, \emph{Finite
subdivision rules},
  Conform. Geom. Dyn. \textbf{5} (2001), 153--196 (electronic).

\bibitem[CFP07]{cfp:subdfromrat}
\bysame, \emph{Constructing subdivision rules from rational maps},
Conform.
  Geom. Dyn. \textbf{11} (2007), 128--136 (electronic).

\bibitem[Cox84]{cox:agm}
David~A. Cox, \emph{The arithmetic-geometric mean of {G}auss},
Enseign. Math.
  (2) \textbf{30} (1984), no.~3-4, 275--330.

\bibitem[CV86]{cullervogtmann:moduli}
Marc Culler and Karen Vogtmann, \emph{Moduli of graphs and
automorphisms of
  free groups}, Invent. Math. \textbf{84} (1986), no.~1, 91--119.

\bibitem[DH84]{DH:orsayI}
Adrien Douady and John~H. Hubbard, \emph{{\'Etude dynamiques des
polyn\^{o}mes
  complexes. (Premi\`ere partie)}}, Publications Mathematiques d'Orsay,
  vol.~02, Universit\'e de Paris-Sud, 1984.

\bibitem[DH85]{DH:orsayII}
\bysame, \emph{{\'Etude dynamiques des polyn\^{o}mes complexes.
(Deuxi\`eme
  partie)}}, Publications Mathematiques d'Orsay, vol.~04, Universit\'e de
  Paris-Sud, 1985.

\bibitem[DH93]{DH:Thurston}
Adrien Douady and John~H. Hubbard, \emph{A proof of {Thurston's}
topological
  characterization of rational functions}, Acta Math. \textbf{171} (1993),
  no.~2, 263--297.

\bibitem[Eil74]{eil}
Samuel Eilenberg, \emph{Automata, languages and machines},
vol.~{A}, Academic
  Press, New York, London, 1974.

\bibitem[Fri87]{fried}
David~L. Fried, \emph{Finitely presented dynamical systems},
Ergod. Th. Dynam.
  Sys. \textbf{7} (1987), 489--507.

\bibitem[FS92]{fornsibon:crfin}
J.~E. Forn{\ae}ss and N.~Sibony, \emph{Critically finite rational
maps on
  {$\mathbb{P}^2$}}, The Madison Symposium on Complex Analysis (Madison, WI,
  1991), Contemp. Math., vol. 137, Amer. Math. Soc., Providence, RI, 1992,
  pp.~245--260.

\bibitem[GM05]{glasnermozes}
Yair Glasner and Shahar Mozes, \emph{Automata and square
complexes}, Geom.
  Dedicata \textbf{111} (2005), 43--64.

\bibitem[GNS00]{grineksu_en}
Rostislav~I. Grigorchuk, Volodymyr~V. Nekrashevich, and
Vitali{\u\i}~I.
  Sushchanskii, \emph{Automata, dynamical systems and groups}, Proceedings of
  the Steklov Institute of Mathematics \textbf{231} (2000), 128--203.

\bibitem[Gro87]{gro:hyperb}
Mikhael Gromov, \emph{Hyperbolic groups}, Essays in Group Theory
(S.~M.
  Gersten, ed.), M.S.R.I. Pub., no.~8, Springer, 1987, pp.~75--263.

\bibitem[IS08]{ishiismillie}
Yutaka Ishii and John Smillie, \emph{Homotopy shadowing},
(preprint
  arXiv:0804.4629v1), 2008.

\bibitem[Jea07]{jeandel:topautom}
Emmanuel Jeandel, \emph{Topological automata}, Theory Comput.
Syst. \textbf{40}
  (2007), no.~4, 397--407.

\bibitem[Kat04]{katsura:generalizing}
Takeshi Katsura, \emph{A class of {$C\sp\ast$}-algebras
generalizing both graph
  algebras and homeomorphism {$C\sp \ast$}-algebras. {I}. {F}undamental
  results}, Trans. Amer. Math. Soc. \textbf{356} (2004), no.~11, 4287--4322
  (electronic).

\bibitem[Kat06a]{katsura:examples}
\bysame, \emph{A class of {$C\sp *$}-algebras generalizing both
graph algebras
  and homeomorphism {$C\sp *$}-algebras. {II}. {E}xamples}, Internat. J. Math.
  \textbf{17} (2006), no.~7, 791--833.

\bibitem[Kat06b]{katsura:ideals}
\bysame, \emph{A class of {$C\sp *$}-algebras generalizing both
graph algebras
  and homeomorphism {$C\sp *$}-algebras. {III}. {I}deal structures}, Ergodic
  Theory Dynam. Systems \textbf{26} (2006), no.~6, 1805--1854.

\bibitem[Kat08]{katsura:pureinfiniteness}
\bysame, \emph{A class of {$C\sp *$}-algebras generalizing both
graph algebras
  and homeomorphism {$C\sp *$}-algebras. {IV}. {P}ure infiniteness}, J. Funct.
  Anal. \textbf{254} (2008), no.~5, 1161--1187.

\bibitem[KL05]{kaimlyubich:laminations}
Vadim~A. Kaimanovich and Mikhail Lyubich, \emph{Conformal and
harmonic measures
  on laminations associated with rational maps}, vol. 173, Memoirs of the
  A.M.S., no. 820, A.M.S., Providence, Rhode Island, 2005.

\bibitem[Koc07]{koch:french}
S.~Koch, \emph{{Teichm\"uller} theory and endomorphisms of
{$\mathbb{P}^n$}},
  Ph.D. thesis, Universit\'e de Provence, 2007.

\bibitem[LM97]{lyubichminsk}
Mikhail Lyubich and Yair Minsky, \emph{Laminations in holomorphic
dynamics}, J.
  Differ. Geom. \textbf{47} (1997), no.~1, 17--94.

\bibitem[LMZ94]{lmz}
Alexander Lubotzky, Shahar Mozes, and Robert~J. Zimmer,
\emph{Superrigidity for
  the commensurability group of tree lattices}, Comment.\ Math.\ Helvetici
  \textbf{69} (1994), 523--548.

\bibitem[MNS00]{mns_en}
Olga Macedo\'nska, Volodymyr~V. Nekrashevych, and Vitali{\u\i}~I.
Sushchansky,
  \emph{Commensurators of groups and reversible automata}, Dopov. Nats. Akad.
  Nauk Ukr., Mat. Pryr. Tekh. Nauky (2000), no.~12, 36--39.

\bibitem[Nek05]{nek:book}
Volodymyr Nekrashevych, \emph{Self-similar groups}, Mathematical
Surveys and
  Monographs, vol. 117, Amer. Math. Soc., Providence, RI, 2005.

\bibitem[Nek08a]{nek:dendrites}
\bysame, \emph{The {Julia} set of a post-critically finite
endomorphism of
  {$\mathbb{CP}^2$}}, (preprint arXiv:0811.2777), 2008.

\bibitem[Nek08b]{nek:filling}
\bysame, \emph{Symbolic dynamics and self-similar groups},
Holomorphic dynamics
  and renormalization. {A volume in honour of John Milnor's 75th birthday}
  (Mikhail Lyubich and Michael Yampolsky, eds.), Fields Institute
  Communications, vol.~53, A.M.S., 2008, pp.~25--73.

\bibitem[Shu69]{shub1}
Michael Shub, \emph{Endomorphisms of compact differentiable
manifolds}, Am. J.
  Math. \textbf{91} (1969), 175--199.

\bibitem[Shu70]{shub2}
\bysame, \emph{Expanding maps}, Global Analysis, Proc. Sympos.
Pure Math.,
  vol.~14, American Math. Soc., Providence, Rhode Island, 1970, pp.~273--276.

\bibitem[VV07]{vorobets:alfree}
Mariya Vorobets and Yaroslav Vorobets, \emph{On a free group of
transformations
  defined by an automaton}, Geometriae Dedicata \textbf{124} (2007), no.~1,
  237--249.

\end{thebibliography}
\end{document}